\documentclass[a4paper,11pt]{article}
\usepackage{mfpaperstuff2}
\usepackage{xfrac}
\usepackage{mathdots}
\usepackage[nottoc,notlot,notlof]{tocbibind}
\usepackage[labelsep=none]{caption}
\begin{document}
\newcommand\isep{$I$-separated\xspace}
\newcommand\rgp[2]{#1\mathbin{\raisebox{-3pt}{\tikz[scale=.3]{\draw(0,0)--++(0,1);\draw(1,0)--++(0,1);\draw(.5,.5)node{\footnotesize$\vphantom{c}\smash{\ol c}$};}}}#2}
\newcommand\rsp[2]{#1\mathbin{\raisebox{-3pt}{\tikz[scale=.3]{\draw(0,0)--++(1,0);\draw(0,1)--++(1,0);\draw(.5,.5)node{\footnotesize$c$};}}}#2}
\newcommand\uiptn{$(u,I)$-partition}
\newcommand\uip{\uiptn\xspace}
\newcommand\uips{\uiptn s\xspace}
\newcommand\ol\widebar
\newcommand\Abdisp{Abacus display}
\newcommand\abdisp{abacus display}
\newcommand\abd{\abdisp\xspace}
\newcommand\abds{\abdisp s\xspace}
\newcommand\Abds{\Abdisp s\xspace}
\newcommand\blo{\operatorname{B}(\La_0)}
\newcommand\sgn{\operatorname{sgn}}
\newcommand\jms[1]{\operatorname{D}^{#1}}
\newcommand\wht{\Yfillcolour{white}}
\newcommand\gry{\Yfillcolour{gray}}
\newcommand\reg{^{\operatorname{reg}}}
\newcommand\tsfrac[2]{\text{\large\sfrac{#1}{#2}}}
\newcommand\mul[1]{\operatorname{m}_{#1}}
\newcommand\mule{\mul e}
\newcommand\Mule{\operatorname{M}_e}
\newcommand\kule{\operatorname{K}_e}
\newcommand\sule{\operatorname{S}_e}
\newcommand\xule{\operatorname{X}_e}
\newcommand\erp{$e$-regular partition\xspace}
\newcommand\hsle{\widehat{\mathfrak{sl}}_e}
\newcommand\dep[1]{\operatorname{dep}(#1)}
\newcommand\add[3]{\operatorname{add}_{#1}^{#2}(#3)}
\newcommand\rem[3]{\operatorname{rem}_{#1}^{#2}(#3)}
\newcommand\yb[1]{\operatorname{bad}(#1)}
\newcommand\eareg{$(e,y)$-regular\xspace}
\newcommand\regn{regularisation\xspace}
\newcommand\drg[2]{$(#1,#2)$-restricted\xspace}
\newcommand\eadereg{\drg ey}
\newcommand\deregn{restrict\-isation\xspace}
\newcommand\earegn{$(e,y)$-\regn}
\newcommand\eaderegn{$(e,y)$-\deregn}
\newcommand\redu[1]{\operatorname{red}(#1)}
\newcommand\zez{\bbz/e\bbz}
\newcommand\ee[1]{\mathbf{e}_{#1}}
\newcommand\ff[1]{\mathbf{f}_{#1}}
\newcommand\eea[1]{\ee{#1}^A}
\newcommand\ffa[1]{\ff{#1}^A}
\newcommand\eeb[1]{\ee{#1}^B}
\newcommand\ffb[1]{\ff{#1}^B}
\newcommand\cl{\prec}
\newcommand\cg{\succ}
\newcommand\cls{\preccurlyeq}
\newcommand\rg[1]{\calr_{#1}}
\newcommand\rgr[2]{\calr_{#1}^{(#2)}}
\newcommand\rega{\rg A}
\newcommand\regb{\rg B}
\newcommand\regar[1]{\rega^{(#1)}}
\newcommand\regbr[1]{\regb^{(#1)}}
\newcommand\regan{\regar n}
\newcommand\regbn{\regbr n}
\newcommand\ladd[1]{\call(#1)}

\title{Crystals, regularisation and the Mullineux map}


\toptitle

\begin{abstract}
The Mullineux map is a combinatorial function on partitions which describes the effect of tensoring a simple module for the symmetric group in characteristic $p$ with the one-dimensional sign representation. It can also be interpreted as an isomorphism between crystal graphs for $\widehat{\mathfrak{sl}}_p$. We give a new combinatorial description of the Mullineux map by expressing this crystal isomorphism as a composition of isomorphisms between different crystals. These isomorphisms are defined in terms of new generalised regularisation maps introduced by Millan Berdasco.

We then given two applications of our new realisation of the Mullineux map, by providing purely combinatorial proofs of a conjecture of Lyle relating the Mullineux map with regularisation, and a theorem of Paget describing the Mullineux map in RoCK blocks of symmetric groups.
\end{abstract}

\tableofcontents

\section{Introduction}

In the representation theory of Kac--Moody algebras and their quantised enveloping algebras, an important role is played by \emph{crystal bases} for integrable modules, and the associated \emph{crystal graphs}, which facilitate a combinatorial approach to studying these modules. An important theme in combinatorial representation theory is to construct a combinatorial model for the crystal of a given module, in which vertices are labelled by simple combinatorial objects, with a combinatorial rule to determine where the arrows go.

In this paper we concentrate on the simplest case of an affine Kac--Moody algebra, namely the algebra $\widehat{\mathfrak{sl}}_e$ for a fixed integer $e\gs2$. This algebra is implicated in the representation theory of Iwahori--Hecke algebras of type~$A$ in quantum characteristic~$e$. In particular, the crystal~$\blo$ of the irreducible highest-weight module with highest weight~$\La_0$ underlies the Brundan--Kleshchev ``modular branching rules'' for these Hecke algebras. This application highlights two particular combinatorial models of~$\blo$: the Misra--Miwa model \cite{mimi}, in which the vertices are labelled by $e$-restricted partitions, and its dual, in which the vertices are labelled by $e$-regular partitions.

In \cite{berg}, Berg found a third model for~$\blo$, called the \emph{ladder crystal}, and showed that it is isomorphic to the $e$-regular model, with an isomorphism being given by James's $e$-regularisation function. In \cite{mfcry}, the author showed that (provided $e\gs3$) these three models are members of an infinite family of models based on partitions. The first aim of the present paper is to provide a new proof of the main result of \cite{mfcry}, by constructing explicit isomorphisms between the crystals in this family. This is done by generalising Berg's proof using a new generalisation of $e$-regularisation due to Millan Berdasco.

We then apply these results to give a new description of the \emph{Mullineux map}. This is an involution on the set of $e$-regular partitions which arises in the modular representation theory of symmetric groups, where it describes the effect of tensoring an irreducible module with the one-dimensional sign module. The \emph{Mullineux problem} asks for a combinatorial description of this map, and several solutions to this problem are now known. The Mullineux map can also be understood in the crystal setting: via the Brundan--Kleshchev branching rules, the Mullineux map describes the isomorphism between the $e$-restricted and $e$-regular models for $\blo$. We construct this isomorphism as a composition of regularisation isomorphisms between the different crystals in our family, thereby giving another solution to the Mullineux problem. We use this new algorithm to give a new (and purely combinatorial) proof of a conjecture of Lyle relating the Mullineux map and regularisation, and a substantial generalisation of a theorem of Paget describing the effect of the Mullineux map in RoCK blocks.

\section{Partitions}

In this section we review some standard background on partitions, before giving some new results due to Millan Berdasco on regularisation.

\subsection{Elementary notation}

We begin with some standard notation. We write~$\bbn$ for the set of positive integers, and~$\bbn_0$ for the set of non-negative integers. Given a subset $\call\subseteq\bbz^2$ and an element $(s,d)\in\bbz^2$, we will write $\call+(s,d)$ to mean $\lset{(r+s,c+d)}{(r,c)\in\call}$.

Throughout this paper, we fix an integer $e\gs2$. We write~$\zez$ for the set of cosets $a+e\bbz$ for $a\in\bbz$. (We do not employ the usual abuse of notation where $\zez$ is identified with the set $\{0,\dots,e-1\}$.)

\subsection{Partitions and Young diagrams}

A \emph{partition} means an infinite weakly decreasing sequence $\la=(\la_1,\la_2,\dots)$ of non-negative integers such that the sum $\card\la=\la_1+\la_2+\dots$ is finite. When writing partitions, we usually group together equal parts and omit trailing zeroes. The partition $(0,0,\dots)$ is written as~$\varnothing$.

We say that a partition~$\la$ is \emph{$e$-regular} if it does not have~$e$ equal positive parts, and that~$\la$ is \emph{$e$-restricted} if $\la_r-\la_{r+1}<e$ for all~$r$.

If~$\la$ is a partition, the \emph{Young diagram} of~$\la$ is the set
\[
\lset{(r,c)\in\bbn^2}{c\ls\la_r}
\]
whose elements we call the \emph{nodes} of~$\la$. In general, a \emph{node} means an element of~$\bbn^2$.

We may abuse notation and identify~$\la$ with its Young diagram; so for example if~$\la$ and~$\mu$ are partitions we may write $\la\subseteq\mu$ to mean that $\la_r\ls\mu_r$ for all~$r$. We draw a Young diagram as an array of boxes in the plane using the English convention, in which the first coordinate increases down the page and the second increases from left to right. We use words such as ``higher'' and ``lower'' with this convention in mind, so that the node $(1,2)$ is higher than, and to the right of, the node $(2,1)$.

A node $(r,c)$ of~$\la$ is \emph{removable} if it can be removed from the Young diagram of~$\la$ to leave the Young diagram of a smaller partition, i.e.\ if $c=\la_r>\la_{r+1}$. A node which does not belong to~$\la$ is an \emph{addable node} of~$\la$ if it can be added to~$\la$ to leave the Young diagram of a larger partition. If $(r,c)$ is an addable node of~$\la$, then we may write $\la\cup(r,c)$ for the partition (whose Young diagram is) obtained by adding this node; similarly we may write $\la\setminus(r,c)$ for the partition obtained by removing a removable node $(r,c)$.

The \emph{residue} of $(r,c)\in\bbz^2$ is defined to be $c-r+e\bbz$. If $i\in\zez$, then we use the term \emph{$i$-node} to mean a node of residue~$i$. We define the \emph{$e$-content} of a partition to be the multiset of the residues of its nodes.

\subsection{Conjugation and the dominance order}

If~$\la$ is a partition, the conjugate partition~$\la'$ is the partition whose Young diagram is obtained by reflecting the Young diagram of~$\la$ in the main diagonal; in other words,~$\la'$ is the partition with $\la'_r=\cardx{\lset{c\in\bbn}{\la_c\gs r}}$ for each~$r$.

The \emph{dominance order} is a natural partial order on partitions: we say that~$\la$ \emph{dominates}~$\mu$, and write $\la\dom\mu$, if $\card\la=\card\mu$ and
\[
\la_1+\dots+\la_r\gs\mu_1+\dots+\mu_r
\]
for all~$r$. Another way to describe this is to say that $\la\dom\mu$ if~$\mu$ can be obtained from~$\la$ by replacing some of the nodes of~$\la$ with lower nodes.

We will use the well-known simple result that conjugation reverses the dominance order: $\la\dom\mu$ \iff $\la'\domby\mu'$.

\subsection{Hooks}

Suppose $(r,c)$ is a node of a partition~$\la$. The \emph{$(r,c)$-hook} of~$\la$ is defined to be the set of all nodes of~$\la$ directly to the right of or directly below $(r,c)$, including $(r,c)$ itself. The \emph{length} of this hook is the number of nodes it contains, i.e.\ $\la_r-c+\la'_c-r+1$, and we refer to the hook as an \emph{$e$-hook} it is has length~$e$. The \emph{arm length} of the hook is the number of nodes to the right of $(r,c)$, i.e.\ $\la_r-c$, and the \emph{leg length} is $\la'_c-r$. The \emph{hand node} of the hook is the node $(r,\la_r)$, and the \emph{foot node} is $(\la'_c,c)$.

For example, the $(2,3)$-hook of the partition $(4^3,3,1)$ has length~$4$, arm length~$1$, leg length~$2$, hand node $(2,4)$ and foot node $(4,3)$, as we see from the following Young diagram.
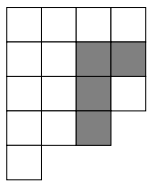
\[
\gyoung(!\wht\ \ \ \ ,\ \ !\gry\ \ ,!\wht\ \ !\gry\ !\wht\ ,\ \ !\gry\ ,!\wht\ )
\]

\emph{Removing} a hook~$H$ from a Young diagram means deleting all the nodes in the hook, and then moving all the nodes below and to the right of~$H$ diagonally up and to the left to create a new Young diagram. (Equivalently, we can remove the nodes of the corresponding \emph{rim hook}, i.e.\ all the nodes along the bottom-right edge of the Young diagram from the hand node of~$H$ to the foot node.) For example, removing the $(2,3)$-hook of the partition $(4^3,3,1)$ from the last example results in the partition $(4,3,2^2,1)$.

\subsection{The abacus}

Combinatorics of partitions is often facilitated by using the \abds introduced by James. We construct an abacus with~$e$ vertical runners labelled with the elements of $\zez$, with the labels $0+e\bbz,1+e\bbz,\dots,-1+e\bbz$ occurring from left to right. We mark positions on these runners labelled with the non-negative integers, reading from left to right along successive rows from the top down. For example, when $e=5$, the positions are labelled as follows.
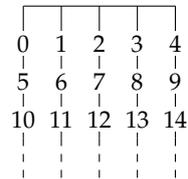
\[
\begin{tikzpicture}[scale=.5,yscale=-1,every node/.style={fill=white,inner sep=1.5pt}]
\draw(0,-1)--++(4,0);
\foreach\x in {0,1,2,3,4}{\draw(\x,-1)--++(0,3.5);\draw(\x,2.5)[dashed]--++(0,1.25);}
\draw(0,0)node{\footnotesize$0$};
\draw(1,0)node{\footnotesize$1$};
\draw(2,0)node{\footnotesize$2$};
\draw(3,0)node{\footnotesize$3$};
\draw(4,0)node{\footnotesize$4$};
\draw(0,1)node{\footnotesize$5$};
\draw(1,1)node{\footnotesize$6$};
\draw(2,1)node{\footnotesize$7$};
\draw(3,1)node{\footnotesize$8$};
\draw(4,1)node{\footnotesize$9$};
\draw(0,2)node{\footnotesize$10$};
\draw(1,2)node{\footnotesize$11$};
\draw(2,2)node{\footnotesize$12$};
\draw(3,2)node{\footnotesize$13$};
\draw(4,2)node{\footnotesize$14$};
\end{tikzpicture}
\]
Given two positions $a<b$ on the abacus, we will say that~$a$ comes \emph{before}, or \emph{earlier than},~$b$, and that~$b$ comes \emph{after}, or \emph{later than},~$a$.

Given a partition~$\la$ and an integer $n\gs\la'_1$, we define the \emph{$n$-bead \abd} for~$\la$ by placing a bead on the abacus at position $\la_r+r-1$ for $r=1,\dots,n$. For example, when $e=5$, the $7$-bead \abd for the partition $(6,4,2,1^2)$ is as follows.
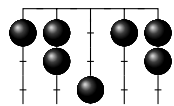
\[
\abacus(lmmmr,bbnbb,nbnnb,nnbnn)
\]
In an \abd, we say that a position is \emph{occupied} if there is a bead at that position, and \emph{empty} otherwise. When we draw \abds, all positions below those drawn should be interpreted as being empty. Every configuration of finitely many beads on the abacus uniquely defines a partition: if the beads are in positions $b_1>\dots>b_n$, then the corresponding partition~$\la$ is given by
\[
\la_r=
\begin{cases}
b_r+r-n&\text{for }r=1,\dots,n
\\
0&\text{ for }r>n.
\end{cases}
\]

An \abd for a partition~$\la$ is particularly useful for visualising removal of hooks: the hooks of length~$e$ in~$\la$ correspond to the positions $b\gs e$ such that position~$b$ is occupied in the \abd while position $b-e$ is empty. Removing the hook corresponds to sliding the bead up from position~$b$ to position $b-e$. The leg length of the hook is the number of occupied positions in the range $[b-e+1,b-1]$, and the arm length is the number of empty positions in $[b-e+1,b-1]$.

It is well-known that \abds behave well with regard to conjugation. Given the $n$-bead \abd for a partition~$\la$, we take a large integer~$m$, and let $b_1,\dots,b_{m-n}$ be the empty positions before position~$m$. Now if we construct a new \abd in which the occupied positions are $m-b_1,\dots,m-b_{m-n}$, this will be an \abd for~$\la'$. In other words, we can construct an \abd for~$\la'$ from an \abd for~$\la$ by truncating at some point after all the beads have appeared, and then replacing each bead with an empty space and each empty space with a bead and rotating through~$180^\circ$.

For example, taking the above \abd for the partition $\la=(6,4,2,1^2)$ and choosing $m=15$, we obtain the following \abd for $\la'=(5,3,2^2,1^2)$.
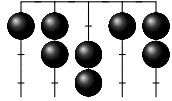
\[
\abacus(lmmmr,bbnbb,nbbnb,nnbnn)
\]

This has the following consequence, which we will need later.

\begin{lemma}\label{ablem}
Suppose~$\la$ is a partition, and take the $n$-bead \abd for~$\la$ with~$e$ runners, for sufficiently large~$n$. Let $t_1<\dots<t_m$ be the first~$m$ empty positions in the \abd, and construct a new \abd by moving the bead at position $t_i-e$ to position $t_i$, for $i=1,\dots,m$ in turn. The resulting configuration is the \abd for the partition obtained by increasing the length of each of the first~$m$ columns of~$\la$ by~$e$.
\end{lemma}

\begin{pf}
Let~$\mu$ be the partition obtained by this procedure. Take an \abd for~$\la$, and let $b_1>\dots>b_m$ be the position of the last~$m$ beads. Then an \abd for~$\mu'$ is obtained by moving the bead at position~$b_i$ to $b_i+e$, for $i=1,\dots,m$ in turn. The definition of the \abd then means that~$\mu'$ is obtained from~$\la'$ by adding~$e$ to each of its first~$m$ parts, which is the same as saying that~$\mu$ is obtained from~$\la$ by adding~$e$ to each of its first~$m$ columns.
\end{pf}

Another important well-known result is the following.

\begin{lemma}\label{contab}
Suppose~$\la$ and~$\mu$ are partitions with the same $e$-content, and construct the $m$-bead \abds for~$\la$ and~$\mu$. Then for each $i\in\zez$ the number of beads on runner~$i$ of the abacus is the same for~$\la$ as it is for~$\mu$.
\end{lemma}

The proof of this \lcnamecref{contab} is just a combination of the fact that two partitions with the same $e$-content have the same $e$-core (a result which goes back to Littlewood \cite{litt}), and the fact that an \abd for the $e$-core of a partition is obtained by sliding all the beads up their runners as far as they will go. This is due to James; we refer the reader to \cite[Section 2.7]{jake} for more information on $e$-cores and the abacus.

\subsection{\earegn}\label{regnsec}

Now we describe some new results due to Millan Berdasco \cite{dmb}, on regularisation of partitions. This generalises $e$-regularisation introduced by James \cite{j1}, and Berg's deregularisation operation \cite{berg}.

Choose an integer $y\in\{1,\dots,e-1\}$. Given $r,c\in\bbz$, we define the \emph{$(e,y)$-ladder}
\[
\ladd{r,c}=\lset{(r+k(y-e),c+ky)}{k\in\bbz}.
\]
The different $(e,y)$-ladders are disjoint, and comprise a partition of~$\bbz^2$. If~$\call$ is an $(e,y)$-ladder and $(s,d)\in\bbz^2$, then $\call+(s,d)$ is also an $(e,y)$-ladder.

We say that two partitions $\la,\mu$ are \emph{$(e,y)$-equivalent} if $\card{\call\cap\la}=\card{\call\cap\mu}$ for every $(e,y)$-ladder~$\call$. This is an equivalence relation on the set of partitions, and we call an equivalence class under this relation an \emph{$(e,y)$-ladder class}.

Now say that a partition is \emph{$(e,y)$-singular} if it has a hook with length~$et$ and arm length $yt-1$ for some $t\gs1$, and \emph{\eareg} otherwise. The next theorem is the main result of \cite{dmb}.

\begin{thm}\label{diego}
Each $(e,y)$-ladder class~$\calc$ contains a unique \eareg partition. This partition dominates every partition in~$\calc$.
\end{thm}

\begin{eg}
The $(3,2)$-ladder class containing the partition $(5,1)$ contains three other partitions. We illustrate these partitions as follows, labelling nodes with letters so that nodes in the same $(3,2)$-ladder get the same label.
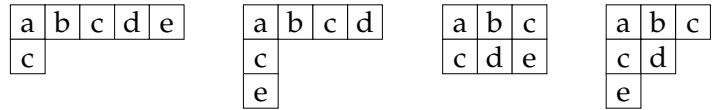
\[\Ystdtext1
\young(abcde,c)\qquad\young(abcd,c,e)\qquad\young(abc,cde)\qquad\young(abc,cd,e)
\]
The unique most dominant partition in this class is $(5,1)$, and this is the only $(3,2)$-regular partition in the class: each of the other partitions has either a $3$-hook with arm length~$1$ or a $6$-hook with arm length~$3$.
\end{eg}

In view of \cref{diego}, we can define the \emph{\earegn} of a partition~$\la$ to be the unique \eareg partition in the same ladder class.

In the case $y=1$, this construction has been known for a long time. The $(e,1)$-regularisation is the same as the $e$-\regn defined by James \cite{j1}, and is constructed by replacing all the nodes of~$\la$ in each ladder with the highest nodes in that ladder. When $y>1$, the construction is not so straightforward: simply replacing the nodes with the highest nodes in their ladders will not in general result in a Young diagram. Berg \cite{berg} addresses the case $y=e-1$, by introducing a notion of ``locked'' nodes of~$\la$, and then constructing the $(e,y)$-\regn of~$\la$ by moving all unlocked nodes to the highest available positions in their ladders. (In fact Berg's convention is the opposite of ours, in that he works with $(e,1)$-ladders and constructs the least dominant partition in the ladder class of~$\la$; but by conjugation this is equivalent to the case $y=e-1$ of \cref{diego} (cf.\ \cref{diegocor} below).)

Millan Berdasco \cite{dmb} gives an algorithm for constructing the $(e,y)$-\regn in general using the abacus, generalising the author's algorithm \cite{mfabreg} for realising the $e$-\regn map on the abacus. This shows in particular that the $(e,y)$-\regn is obtained by moving nodes up their $(e,y)$-ladders. We will see this algorithm in the next section.

We will also need to consider the least dominant partition in each $(e,y)$-ladder class, and we do this by considering conjugate partitions. Define $S'=\lset{(c,r)}{(r,c)\in S}$ for any set $S\subseteq\bbz^2$, and observe that if~$\call$ is an $(e,y)$-ladder, then~$\call'$ is an $(e,e-y)$-ladder. Hence two partitions $\la,\mu$ are $(e,y)$-equivalent \iff~$\la'$ and~$\mu'$ are $(e,e-y)$-equivalent. Furthermore, if~$H$ is a hook with length~$et$ and arm length $yt-1$, then~$H'$ is a hook with length~$et$ and arm length $(e-y)t$. Now say that a partition is \emph{\eadereg} if it has no hooks of length~$et$ and arm length~$yt$ for any~$t$. This generalises the definition of an $e$-restricted partition: ``$e$-restricted'' is the same as ``\drg e{e-1}''.

The following result follows from \cref{diego} and the fact that conjugation reverses the dominance order.

\begin{cory}\label{diegocor}
Each $(e,y)$-ladder class~$\calc$ contains a unique \eadereg partition. This partition is dominated by every partition in~$\calc$.
\end{cory}

The case $y=1$ of this \lcnamecref{diegocor} is due to Berg. If~$\la$ is a partition, then we refer to the least domi\-nant partition in the same $(e,y)$-ladder class as~$\la$ as the \emph{\eaderegn} of~$\la$. (We prefer the awkward artificial word ``restrictisation'' over the more natural ``restriction'', because the latter word is already widely used.)

\begin{eg}
Referring back to the last example, we see that the unique least dominant partition in the $(3,2)$-class containing $(5,1)$ is $(3,2,1)$. This partition is \drg32, while each of the other partitions in the class has a $3$-hook with arm length~$2$.
\end{eg}

\subsection{\earegn on the abacus}\label{regnabsec}

Now we give an algorithm due to Millan Berdasco for computing the \earegn of a partition using the abacus. We will use this in \cref{splitsec} when we consider the Mullineux map on the abacus.

Take an integer $y\in[1,e-1]$ and a partition~$\nu$ which is $(e,y)$-singular. Our task is to construct a more dominant partition which is $(e,y)$-equivalent to~$\mu$; by doing this repeatedly, we will eventually reach the \earegn of~$\nu$.

Construct an \abd for~$\nu$. Because~$\nu$ is $(e,y)$-singular, there is some $t\in\bbn$ such that~$\nu$ has a hook of length~$et$ and arm length $yt-1$. We will assume that the only such~$t$ occurring is $t=1$; if this is not the case, then taking the largest~$t$ that occurs, we can replace $e,y$ with $et,yt$, and apply the same procedure to get a partition $\mu\doms\nu$ which is $(et,yt)$-equivalent to~$\nu$. Because each $(e,y)$-ladder is a union of $(et,yt)$-ladders,~$\mu$ will also be $(e,y)$-equivalent to~$\nu$, as required.

Assuming $t=1$ is the only value of~$t$ that occurs, there is at least one occupied position~$b$ in the abacus such that position $b-e$ is empty and there are exactly~$y$ empty positions in the range $[b-e,b]$. We take the largest such~$b$, and let~$E$ be the union of the congruence classes modulo~$e$ of the empty positions in $[b-e,b]$. Now let $b_1<\dots<b_m$ be the occupied positions in~$E$ after position $b-e$; in particular, $b_1=b$. Also, let $t_1<t_2<\dots$ be the empty positions in $\bbn_0\setminus E$ after position~$b$. Let $d\in\{1,\dots,m\}$ be minimal such that either $t_d<b_{d+1}$ or $d=m$. Now construct a new \abd by moving the bead at position~$b_i$ to position $b_i-e$ and moving the bead at position $t_i-e$ to position~$t_i$ for $i=1,\dots,d$ in turn. This gives the \abd for a partition~$\kappa$, and Millan Berdasco \cite[Proposition 4.1]{dmb} proves that $\kappa\doms\nu$ and~$\kappa$ is $(e,y)$-equivalent to~$\nu$, as required.

\begin{eg}
Take $(e,y)=(5,3)$ and $\la=(9,3^3,2)$, with the following $5$-bead \abd.
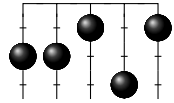
\[
\abacus(lmmmr,nnbnb,bbnnn,nnnbn)
\]
We can see that $b=5$, with $E=(0+5\bbz)\cup(1+5\bbz)\cup(3+5\bbz)$. Hence $(b_1,\dots,b_m)=(5,6,13)$ and $(t_1,t_2,\dots)=(7,9,12,14,\dots)$. Hence $d=2$. So~$\kappa$ is defined by moving the beads at positions~$5$ and~$6$ up and the beads at positions~$7$ and~$9$ down, giving $\kappa=(9,6,5)$.
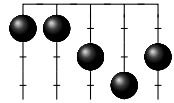
\[
\abacus(lmmmr,bbnnn,nnbnb,nnnbn)
\]
\end{eg}

\section{Crystals}\label{crysec}

Now we introduce crystals, giving a simplified definition suitable for our purposes. We recall the definition of a family of crystals from \cite{mfcry}, and use \earegn to prove that they are isomorphic. Excellent references for crystals are the books by Kashiwara \cite{kashbook}, Hong and Kang \cite{hongkang}, and Bump and Schilling \cite{bs}. Here we provide an abbreviated account, restricting to the special case we need for this paper.

\subsection{Introduction to crystals}

A \emph{$\zez$-crystal} means a set~$B$ (not including~$0$ as an element) together with functions $\ee i,\ff i:B\to B\cup\{0\}$ for each $i\in\zez$ with the property that if $b,c\in B$ and $i\in\zez$ then $\ee ib=c$ \iff $\ff ic=b$. The associated \emph{crystal graph} is a labelled directed graph with~$B$ as its vertex set, and an arrow $b\stackrel i\rightarrow c$ whenever $\ff ib=c$. Often we abuse notation by not distinguishing between a crystal and the crystal graph.

Now let~$\hsle$ denote the affine Kac--Moody algebra of type~$A_{e-1}^{(1)}$ (see \cite[Chapter 2]{hongkang} for an uncomplicated account of the definition and classification of Kac--Moody algebras). An $\hsle$-crystal means a $\zez$-crystal~$B$ endowed with a weight function~$\operatorname{wt}$ from~$B$ to the weight space for~$\hsle$, and functions $\epsilon_i,\phi_i:B\to\bbz\cup\{-\infty\}$ for each $i\in\zez$ which satisfy certain axioms.

Certain types of module for the quantum group~$U_q(\hsle)$ come equipped with crystals defined in a natural way from a crystal basis for the module. We say that an abstract crystal is \emph{regular} (the term \emph{normal} is used in \cite{bs}) if it arises as the crystal of a module in this way. An important theme in the theory of crystals is finding combinatorial models for regular crystals, i.e.\ finding an abstract crystal isomorphic to the crystal of a given module, with simple combinatorial objects (such as partitions or tableaux)  as vertices and a combinatorial rule to determine where the arrows go.

In this paper we are concerned with a family of models for one particular regular crystal: this is the \emph{basic crystal}~$\blo$, which is the crystal of the irreducible integrable highest-weight module~$V(\La_0)$ for~$U_q(\hsle)$. There are two very well known models for this crystal:
\begin{itemize}
\item
the \emph{Misra--Miwa model} \cite{mimi}, whose vertex set is the set of all $e$-restricted partitions;
\item
the dual of the Misra--Miwa model, whose vertex set is the set of all $e$-regular partitions.
\end{itemize}

We will now we assume $e\gs3$; this assumption will be in force until the end of \cref{crysec}, where we will make some comments on the case $e=2$. In the case $e\gs3$, Berg \cite{berg} found another model for~$\blo$ (the ``ladder crystal'') with a different set of partitions as its vertex set, and then in \cite{mfcry} the author showed that all three models are cases of a continuous family of models for~$\blo$. The proof in \cite{mfcry} is long and technical, using results of Stembridge \cite{stem} and Danilov--Karzanov--Koshevoy \cite{dkk} to show that each of the proposed models is a regular crystal by analysing its local structure, and then appealing to the uniqueness of regular crystals with a unique highest-weight vertex.

The object in this paper is to give a more direct proof that the crystals from \cite{mfcry} are models for~$\blo$, by showing directly that they are all isomorphic to the dual of the Misra--Miwa model; this generalises Berg's approach. A by-product of this is a direct construction of the (unique) isomorphism from the Misra--Miwa model to its dual; this provides a new realisation of the \emph{Mullineux map}, which is the subject of \cref{mullsec}.

For this paper, it suffices to think of a crystal simply as a directed graph whose arrows are labelled with elements of~$\zez$. An \emph{isomorphism} between two crystals then just means a bijection~$\alpha$ between their vertex sets such that there is an arrow $b\stackrel i\rightarrow c$ \iff there is an arrow $\alpha(b)\stackrel i\rightarrow\alpha(c)$. The additional functions~$\operatorname{wt},\epsilon_i,\phi_i$ and their compatibility with isomorphisms follow naturally in the cases we are concerned with, so we can ignore them in this paper.

\subsection{A family of crystals}\label{cryfamilysec}

We now define the family of crystals introduced in \cite[Section 2.3]{mfcry}; the Misra--Miwa model and its dual arise as special cases.

Define an \emph{arm sequence} to be a sequence $A=(A_1,A_2,\dots)$ of integers satisfying
\begin{itemize}
\item
$t-1\ls A_t\ls(e-1)t$ for all $t\gs1$, and
\item
$A_{t+u}\in\{A_t+A_u,A_t+A_u+1\}$ for all $t,u\gs1$.
\end{itemize}

In fact arm sequences are easily classified. Given any real number $y\in[1,e-1]$, we define two arm sequences~$A^{y+}$ and~$A^{y-}$ by
\[
A^{y+}_t=\inp{yt},\qquad A^{y-}=\roundup{yt-1}
\]
for $t\gs1$. By \cite[Lemma 7.4]{mfcry} these sequences are arm sequences, and every arm sequence has one of these forms. Obviously if~$y$ is irrational then $A^{y+}=A^{y-}$, but otherwise these arm sequences are distinct.

Now given two arm sequences~$A$ and~$B$, write $A\ls B$ if $A_t\ls B_t$ for all~$t$. It follows from the classification above that~$\ls$ is a total order on the set of arm sequences. Specifically, if $x,y\in[1,e-1]$ with $x<y$, we get $A^{x\pm}<A^{y\pm}$, while $A^{y-}<A^{y+}$ for every rational $y\in[1,e-1]$.

We will define a crystal~$\rega$ for each arm sequence~$A$. First we define the underlying set of partitions which will be the vertex set of~$\rega$. Say that a partition~$\la$ is \emph{$A$-regular} if it has no hook with length~$et$ and arm length~$A_t$ for any~$t$. We define~$\rega$ to be the set of $A$-regular partitions.

\begin{eg}
Suppose $e=4$, and let $\la=(5,2,1^2)$. The hooks of~$\la$ of length divisible by~$4$ are a $4$-hook with arm length~$1$, and an $8$-hook with arm length~$4$. So~$\la$ is $A$-regular for any arm sequence~$A$ beginning $(0,\dots)$ or $(2,5,\dots)$ or $(3,\dots)$.
\end{eg}

Now we need to define the crystal operators~$\ee i$ and~$\ff i$ on~$\rega$, for $i\in\zez$. Recall that an \emph{$i$-node} means a node of residue~$i$. We define a total order (depending on~$A$) on the set of all $i$-nodes. Given two different $i$-nodes $(r,c)$ and $(s,d)$, the integer $s-r+c-d$ equals~$et$ for some integer~$t$. By interchanging the two nodes if necessary, we assume $t\gs0$. Now we set $(r,c)\cl(s,d)$ if $c-d\ls A_t$, and $(r,c)\cg(s,d)$ otherwise. For this purpose, we read~$A_0$ as~$0$. Then it is easy to check that~$\cls$ is a total order on the set of $i$-nodes.

Now take $\la\in\rega$. Let $(r_1,c_1),\dots,(r_m,c_m)$ be the addable and removable $i$-nodes of~$\la$, ordered so that $(r_1,c_1)\cg\cdots\cg(r_m,c_m)$. Define the \emph{$i$-signature} of~$\la$ (with respect to~$A$) to be the sequence of signs $s=s_1\dots s_m$ defined by $s_k=+$ if $(r_k,c_k)$ is an addable node, and $s_k=-$ if $(r_k,c_k)$ is a removable node. The \emph{reduction} of~$s$ is the sequence~$\redu s$ obtained from~$s$ by repeatedly deleting adjacent pairs~$+-$. If there are any~$-$ signs in~$\redu s$, then the removable node $(r_k,c_k)$ corresponding to the last one is the \emph{good} $i$-node of~$\la$, and we define $\ee i\la=\la\setminus(r_k,c_k)$; otherwise we set $\ee i\la=0$. If there are any~$+$ signs in~$\redu s$, then the addable node $(r_l,c_l)$ corresponding to the first one is called the \emph{cogood} $i$-node of~$\la$, and we set $\ff i\la=\la\cup(r_l,c_l)$; otherwise we set $\ff i\la=0$.

It is an easy consequence of the construction that if $\la,\mu\in\rega$ then $\la=\ee i\mu$ \iff $\mu=\ff i\la$. Furthermore, it is shown in \cite[Proposition 6.1]{mfcry} that if $\la\in\rega$, then $\ee i\la,\ff i\la\in\rega\cup\{0\}$. So the functions $\ee i,\ff i$ for $i\in\zez$ endow~$\rega$ with the structure of a crystal. Furthermore, the only source in this crystal (i.e.\ the only vertex with no incoming arrows) is~$\varnothing$, by \cite[Proposition 6.3]{mfcry}.

\begin{eg}
Suppose $e=4$ and~$A$ is an arm sequence. Suppose $\la=(5,2,1^2)$ is $A$-regular. Then $A_1\neq1$ and $A_2\neq4$. Take $i=0+4\bbz$. The addable and removable $i$-nodes of~$\la$ are $(1,5)$, $(2,2)$, and $(5,1)$. The ordering of these nodes, the $i$-signature~$s$, the reduced $i$-signature~$\redu s$ depend on~$A_1$ as follows.
\[
\begin{array}{cccccc}\hline
A_1&\text{order}&s&\redu s&\ee i\la&\ff i\la\\\hline
0&(1,5)\cg(2,2)\cg(5,1)&--+&--+&(5,1^3)&(5,2,1^3)\\
2&(5,1)\cg(1,5)\cg(2,2)&+--&-&(5,1^3)&0\\
3&(5,1)\cg(2,2)\cg(1,5)&+--&-&(4,2,1^2)&0\\\hline
\end{array}
\]
\end{eg}

We now consider the two extreme special cases of this construction. In the special case where $A=A^{(e-1)+}=(e-1,2e-2,3e-3,\dots)$, a partition~$\la$ is $A$-regular \iff it is $e$-restricted: the $e$-restricted condition simply says that~$\la$ has no $e$-hook with arm length $e-1$, but this automatically implies that~$\la$ has no $et$-hook with arm length $(e-1)t$ for any~$t$. If~$\la$ is $e$-restricted, then the ordering $(r_1,c_1)\cg\cdots\cg(r_m,c_m)$ on the addable and removable $i$-nodes is simply the order of these nodes from bottom to top in the Young diagram. This means that $\rg{A^{(e-1)+}}$ is the Misra--Miwa model. Similarly, if $A=A^{1-}=(0,1,2,\dots)$, then a partition is $A$-regular \iff it is $e$-regular, in which case the ordering on the addable and removable $i$-nodes is the order from top to bottom, which means that $\rg{A^{1-}}$ is the dual of the Misra--Miwa model.

\subsection{\earegn{} and crystal isomorphisms}\label{isosec}

Now we come to our first main result, which uses the regularisation and restrictisation operations from \cref{regnsec} to give crystal isomorphisms.

\begin{thm}\label{main}
Suppose $y\in[1,e-1]\cap\bbq$, and let~$z$ be the denominator of~$y$. Then $(ez,yz)$-\regn defines an isomorphism of crystals $\rg{A^{y+}}\to\rg{A^{y-}}$. The inverse isomorphism is given by $(ez,yz)$-\deregn.
\end{thm}

For the remainder of \cref{isosec} we fix $y\in[1,e-1]\cap\bbq$, we let~$z$ denote the denominator of~$y$, and we write $A=A^{y+}$ and $B=A^{y-}$. We will say ``ladder'' to mean ``$(ez,yz)$-ladder'', and use the terms ``ladder equivalent'' and ``ladder class'' similarly.

Proving \cref{main} involves two parts: first showing that $(ez,yz)$-\regn and $(ez,yz)$-\deregn give mutually inverse functions between the sets~$\rega$ and~$\regb$, and then showing that these functions commute with the crystal operators~$\ee i$ and~$\ff i$.

First we show that $(ez,yz)$-\regn maps~$\rega$ to~$\regb$. In the case $z=1$ this is immediate: if~$y$ is an integer, then ``$(e,y)$-regular'' and ``$A^{y-}$-regular'' mean the same thing. But for $z>1$ more work is needed. Let's assume that $z>1$, and say that a hook is \emph{$y$-bad} if it has length~$te$ and arm length $\inp{yt}$ for some integer~$t$ not divisible by~$z$. Let~$\yb\la$ denote the number of $y$-bad hooks of~$\la$.

\begin{propn}\label{samehooks}
Suppose~$\la$ and~$\mu$ lie in the same ladder class. Then $\yb\la=\yb\mu$.
\end{propn}

\begin{pf}
Recall that if $r,c\in\bbz$ then the ladder~$\ladd{r,c}$ is defined as
\[
\ladd{r,c}=\lset{(r+k(y-e)z,c+kyz)}{k\in\bbz}.
\]
We need to show that~$\yb\la$ depends only on the number of nodes of~$\la$ in each ladder, and we do this by induction on~$\card\la$. For the inductive step, take a partition~$\nu$ with an addable node $(r,c)$; then we need to show that $\yb{\nu\cup(r,c)}$ depends only on~$\yb\nu$, the ladder~$\ladd{r,c}$, and the number of nodes of~$\nu$ in each ladder.

Define the set
\[
\calm=\bigcup_{t=1}^{z-1}\ladd{r+A_t-te,c+A_t}.
\]
Observe that because~$z$ is the denominator of~$y$,
\[
A_{t+z}=\inp{y(t+z)}=\inp{yt}+yz=A_t+yz
\]
for any~$t$, and so
\[
\ladd{r+A_{t+z}-(t+z)e,c+A_{t+z}}=\ladd{r+A_t-te,c+A_t}.
\]
So in fact~$\calm$ contains the ladder $\ladd{r+A_t-te,c+A_t}$ for every $t\in\bbn$ not divisible by~$z$.

Now we compare the hooks of $\nu\cup(r,c)$ and~$\nu$. For every hook of~$\nu$ with foot node $(r-1,c)$, there is a hook of $\nu\cup(r,c)$ with the same hand node and with foot node $(r,c)$. Similarly, for every hook of~$\nu$ with hand node $(r,c-1)$, there is a hook of $\nu\cup(r,c)$ with the same foot node and with hand node $(r,c)$. Additionally, $\nu\cup(r,c)$ has a $1$-hook with hand and foot node both equal to $(r,c)$. Otherwise, the hooks of~$\nu$ and $\nu\cup(r,c)$ coincide.

Now observe that:
\begin{itemize}
\item
a hook with foot node $(r,c)$ and hand node $(s,d)$ is $y$-bad \iff $(s-1,d)\in\calm$;
\item
a hook with foot node $(r-1,c)$ and hand node $(s,d)$ is $y$-bad \iff $(s,d)\in\calm$;
\item
a hook with hand node $(r,c)$ and foot node $(s,d)$ is $y$-bad \iff $(s,d-1)\in\calm$;
\item
a hook with hand node $(r,c-1)$ and foot node $(s,d)$ is $y$-bad \iff $(s,d)\in\calm$.
\end{itemize}
So $\yb{\nu\cup(r,c)}-\yb\nu$ equals
\begin{align*}
&\cardx{\lset{(s,d)\in\nu\vphantom{l^2_2}}{s<r,\ (s,d+1)\notin\nu\text{ and }(s-1,d)\in\calm}}\\
-&\cardx{\lset{(s,d)\in\nu\vphantom{l^2_2}}{s<r,\ (s,d+1)\notin\nu\text{ and }(s,d)\in\calm}}\\
+&\cardx{\lset{(s,d)\in\nu\vphantom{l^2_2}}{s>r,\ (s+1,d)\notin\nu\text{ and }(s,d-1)\in\calm}}\\
-&\cardx{\lset{(s,d)\in\nu\vphantom{l^2_2}}{s>r,\ (s+1,d)\notin\nu\text{ and }(s,d)\in\calm}}.
\end{align*}
Taking each $(s,d)\in\calm$ in turn and examining the possible intersections of~$\nu$ with the set $\{(s,d),(s+1,d),(s,d+1),(s+1,d+1)\}$ and its contribution to the above sum, we find that $\yb{\nu\cup(r,c)}-\yb\nu$ is the number of configurations
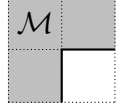
\[
\begin{tikzpicture}[scale=.7,label distance=0pt]
\draw[gray!50!white,fill=white!50!gray](0,0)--++(1,0)--++(0,1)--++(1,0)--++(0,1)--++(-2,0)--++(0,-2);
\foreach\x in {0,1,2}\draw[densely dotted](0,\x)--++(2,0)(\x,0)--++(0,2);
\draw[thick](1,0)--++(0,1)--++(1,0);
\draw(.5,1.5)node{$\calm$};
\end{tikzpicture}
\]
in~$\nu$ minus the number of configurations
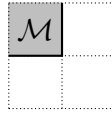
\[
\begin{tikzpicture}[scale=.7,label distance=0pt]
\draw[gray!50!white,fill=white!50!gray](0,1)--++(1,0)--++(0,1)--++(-1,0)--++(0,-1);
\foreach\x in {0,1,2}\draw[densely dotted](0,\x)--++(2,0)(\x,0)--++(0,2);
\draw[thick](0,1)--++(1,0)--++(0,1);
\draw(.5,1.5)node{$\calm$};
\end{tikzpicture}
\]
in~$\nu$. In these diagrams, the box marked~$\calm$ indicates an element of~$\calm$ (either higher or lower than the node $(r,c)$), the shaded boxes are nodes of~$\nu$, and the unshaded boxes are nodes not lying in~$\nu$. (We include the case of boxes $(s,d)\in\calm$ with~$s$ or~$d$ equal to~$0$; for these cases, the shaded boxes can include elements of $\bbz^2\setminus\bbn^2$.)

Now we obtain
\[
\yb{\nu\cup(r,c)}=\yb\nu+\card{\nu\cap\calm}-\card{\nu\cap(\calm+(0,1))}-\card{\nu\cap(\calm+(1,0))}+\card{\nu\cap(\calm+(1,1))}.
\]
Since each of the sets~$\calm$, $\calm+(0,1)$, $\calm+(1,0)$, $\calm+(1,1)$ is a union of ladders, this is all we need for our inductive step.
\end{pf}

\begin{cory}\label{bijec}
$(ez,yz)$-\regn yields a bijection $\rega\to\calr_B$, with inverse given by $(ez,yz)$-\deregn.
\end{cory}

\begin{pf}
Take $\la\in\rega$. Then~$\la$ has no hooks of length~$et$ and arm length $A_t=\inp{yt}$ for any~$t$, and so in particular has no $y$-bad hooks. Now let~$\mu$ be the $(ez,yz)$-\regn of~$\la$. Then by \cref{diego}~$\mu$ has no hooks of length~$ezt$ and arm length $yzt-1$ for any~$t$, and by \cref{samehooks}~$\mu$ has no $y$-bad hooks. Hence~$\mu$ has no hooks of length~$et$ and arm length~$B_t$ for any~$t$, i.e.\ $\mu\in\calr_B$.

So $(ez,yz)$-\regn gives a function $\rega\to\regb$. By conjugating everything, we can show in the same way that $(ez,yz)$-\deregn gives a function from~$\calr_B$ to~$\rega$. Every partition in~$\regb$ is $(ez,yz)$-regular, so equals the $(ez,yz)$-\regn of its $(ez,yz)$-\deregn. Similarly, every partition in~$\rega$ equals the $(ez,yz)$-\deregn of its $(ez,yz)$-\regn. So the two functions $\rega\leftrightarrow\regb$ are inverses of each other.
\end{pf}

Now we come to the second part of the proof of \cref{main}: showing that the bijections between~$\rega$ and~$\regb$ preserve the crystal operators~$\ee i$ and~$\ff i$ for each~$i$. For clarity, we will write~$\eea i$ and~$\ffa i$ for the crystal operators on~$\rega$, and define~$\eeb i$ and~$\ffb i$ similarly.

We continue to use ``ladder'' to mean ``$(ez,yz)$-ladder''. Clearly all the nodes in a given ladder have the same residue, and we say that a ladder in which all nodes have residue~$i$ is an $i$-ladder.

If $(r,c)\in\bbz^2$, then we define the \emph{depth} $\dep{r,c}=yr+(e-y)c$.

\begin{lemmac}{dmb}{Lemma 3.2}\label{laddepth}
Suppose $(r,c),(s,d)\in\bbz^2$. Then $(r,c)$ and $(s,d)$ lie in the same ladder \iff they have the same depth and the same residue.
\end{lemmac}

\begin{pf}
It is trivial to check that if $(r,c)$ and $(s,d)$ lie in the same ladder then they have the same depth and the same residue. Conversely, suppose $(r,c)$ and $(s,d)$ have the same residue and the same depth. The latter statement says $yr+(e-y)c=ys+(e-y)d$, which rearranges to
\[
yz\frac{r-c-s+d}e=z(d-c).
\]
The factor $(r-c-s+d)/e$ is an integer because $(r,c)$ and $(s,d)$ have the same residue. By definition~$z$ and~$yz$ are coprime integers, so~$yz$ divides~$d-c$; let's say $d-c=kyz$ for $k\in\bbz$. But then $(s,d)=(r,c)+k(yz-ez,yz)$, so $(r,c)$ and $(s,d)$ lie in the same ladder.
\end{pf}

(In fact, the depth condition in \cref{laddepth} is used in \cite{dmb} as the definition of a ladder.) \cref{laddepth} allows us to impose a total order on the set of $i$-ladders: we set $\calk<\call$  if the depth of the nodes in~$\calk$ is less than the depth of the nodes in~$\call$. In particular, given an $i$-ladder~$\call$, the $i$-ladder $\call+(1,1)$ satisfies $\call<\call+(1,1)$.

Now we study the relationship between ladders and $i$-signatures. First we need a lemma.

\begin{lemma}\label{noaddrem}
Suppose~$\la$ is a partition with a removable $i$-node $(r,c)$ and an addable $i$-node $(s,d)$ with $\dep{s,d}-e<\dep{r,c}<\dep{s,d}$. Then~$\la$ has a $y$-bad hook.
\end{lemma}

\begin{pf}
Let's assume that $c<d$ (the case $c>d$ is very similar). The fact that $\dep{s,d}-e<\dep{r,c}$ implies that $s<r$, and we claim that the hook with hand node $(s,d-1)$ and foot node $(r,c)$ is $y$-bad. The arm length of this hook is $d-c-1$, and its length is $d-c-s+r$, which is divisible by~$e$ because $(r,c)$ and $(s,d)$ are both $i$-nodes.

The assumption that $\dep{s,d}-e<\dep{r,c}<\dep{s,d}$ says
\[
ys+(e-y)d-e<yr+(e-y)c<ys+(e-y)d
\]
which rearranges to
\[
d-c-1<\frac{y(d-c-s+r)}e<d-c.
\]
This says in particular that $y(d-c-s+r)/e$ is not an integer, so that the integer $t=(d-c-s+r)/e$ is not divisible by~$z$. Furthermore, $d-c-1=\inp{yt}=A_t$, so the hook is~$y$-bad.
\end{pf}

Now for any partition~$\la$ and any $u\in\bbq$, write~$\rem iu\la$ for the number of removable $i$-nodes of~$\la$ of depth~$u$, and~$\add iu\la$ for the number of addable $i$-nodes of~$\la$ of depth~$u$. We want to compare these numbers for two partitions which are ladder-equivalent.

\begin{lemma}\label{addrem}
Suppose~$\la$ and~$\mu$ are two ladder-equivalent partitions, and that $i\in\zez$ and $u\in\bbq$. Then
\[
\add iu\la-\rem i{u-e}\la=\add iu\mu-\rem i{u-e}\mu.
\]
\end{lemma}

\begin{pf}
We prove this \lcnamecref{addrem} by induction on~$\card\la$. If there are no $i$-nodes of depth~$u$ then the result is trivial, so assume there are such nodes, and let~$\call$ be the ladder containing them. Then the ladder containing the $i$-nodes of depth $u-e$ (if there are any) is $\call-(1,1)$. To prove our result by induction, it suffices to take a partition~$\nu$ with an addable node $(r,c)$, and show that the difference
\[
D=\add iu{\nu\cup(r,c)}-\rem i{u-e}{\nu\cup(r,c)}-\add iu\nu+\rem i{u-e}\nu
\]
depends only on the ladder that contains $(r,c)$. Adding the node $(r,c)$ to~$\nu$ can only affect the addable nodes in~$\call$ if $(r,c)$ lies in~$\call$, $\call-(0,1)$ or $\call-(1,0)$. Similarly, adding $(r,c)$ can only affect the removable nodes in $\call-(1,1)$ if $(r,c)$ lies in $\call-(0,1)$, $\call-(1,0)$ or $\call-(1,1)$. So we have four cases to consider.
\begin{description}
\item[$(r,c)\in\call$]\indent\\In this case $\add iu{\nu\cup(r,c)}=\add iu\nu-1$, because $(r,c)$ is an addable node of~$\nu$ but not of $\nu\cup(r,c)$. So $D=-1$.
\item[$(r,c)\in\call-(1,0)$]\indent\\In this case consider the node $(r+1,c-1)$.
\begin{itemize}
\item
If $(r+1,c-1)\in\nu$, then $\add iu{\nu\cup(r,c)}=\add iu\nu+1$, because the node $(r+1,c)$ is an addable node of $\nu\cup(r,c)$ but not of~$\nu$; but $\rem i{u-e}{\nu\cup(r,c)}=\rem i{u-e}\nu$.
\item
If $(r+1,c-1)\notin\nu$, then $\add iu{\nu\cup(r,c)}=\add iu\nu$ but $\rem i{u-e}{\nu\cup(r,c)}=\rem i{u-e}\nu-1$, because $(r,c-1)$ is a removable node of~$\nu$ but not of $\nu\cup(r,c)$.
\end{itemize}
Either way, $D=1$ in this case.
\item[$(r,c)\in\call-(0,1)$]\indent\\Similarly to the last case, we always have $D=1$ in this case.
\item[$(r,c)\in\call-(1,1)$]\indent\\In this case $D=-1$, similarly to the first case.\qedhere
\end{description}
\end{pf}

Now we can consider $i$-signatures. Recall from \cref{cryfamilysec} that the definition of the $i$-signature of a partition~$\la$ depends on the order~$\cls$ on the set of $i$-nodes, which in turn depends on the arm sequence~$A$. Since we are concerned with two different arm sequences $A=A^{y+}$ and $B=A^{y-}$, we will write~$\cls_A$ and~$\cls_B$ for the associated orders, and refer to the $(A,i)$-signature and the $(B,i)$-signature of a partition, and to $A$-good and $B$-good $i$-nodes.
Here is the crucial result.

\begin{propn}\label{samredsig}
Suppose $\la\in\rega$, and let~$\mu$ be the $(ez,yz)$-\regn of~$\la$.
\begin{enumerate}
\item
The reduced $(A,i)$-signature of~$\la$ equals the reduced $(B,i)$-signature of~$\mu$.
\item
If~$\la$ has an $A$-good $i$-node, then~$\mu$ has a $B$-good $i$-node in the same ladder.
\end{enumerate}
\end{propn}

\begin{pf}
First suppose $(r,c)$ and $(s,d)$ are $i$-nodes. If $\dep{r,c}<\dep{s,d}$, then $(r,c)\cl_A(s,d)$. So when we calculate the $(A,i)$-signature of~$\la$, the signs~$\pm$ corresponding to ladders of greater depth come before those of smaller depth.

If $\dep{r,c}=\dep{s,d}$, then $(r,c)\cl_A(s,d)$ \iff $r<s$. Now we observe that if~$\la$ has an addable $i$-node $(r,c)$ and a removable $i$-node $(s,d)$ of the same depth, then $r<s$ (otherwise the hook with foot node $(s,d)$ and hand node $(r-1,c)$ would have length~$et$ and arm length~$A_t$ for some~$t$, contradicting the assumption that $\la\in\rega$). What this means is that in the $(A,i)$-signature of~$\la$, the~$-$ signs corresponding to removable nodes in a given ladder precede the~$+$ signs corresponding to addable nodes in the same ladder.

So we can write down the $(A,i)$-signature of~$\la$, in terms of the addable and removable nodes in each ladder. Given an integer~$n$, we let~$[n]$ denote a string of~$+$ signs of length~$n$ if $n\gs0$, or a string of~$-$ signs of length~$-n$ if $n<0$. Now if we let $u_1<u_2<\dots$ be the possible depths that $i$-nodes can have, then the first part of the proof shows that the $(A,i)$-signature of~$\la$ equals the concatenation
\[
\dots\big[{-}\rem i{u_3}\la\big]\big[\add i{u_3}\la\big]\big[{-}\rem i{u_2}\la\big]\big[\add i{u_2}\la\big]\big[{-}\rem i{u_1}\la\big]\big[\add i{u_1}\la\big].
\]
Because~$\la$ has no $y$-bad hooks, \cref{noaddrem} shows that if $v-e<u<v$ then at least one of~$\rem iu\la$ and~$\add iv\la$ is zero, so the string $[\add iv\la][-\rem iu\la]$ is the same as the string $[-\rem iu\la][\add iv\la]$. Applying this observation repeatedly in the above expression for the $(A,i)$-signature of~$\la$, we find that this signature equals
\[
\dots\big[\add i{u_3}\la\big]\big[{-}\rem i{u_3-e}\la\big]\big[\add i{u_2}\la\big]\big[{-}\rem i{u_2-e}\la\big]\big[\add i{u_1}\la\big]\big[{-}\rem i{u_1-e}\la\big].
\]
(Here we have introduced some empty strings $\rem i{u-e}\la$ for those~$u$ such that there are nodes of depth~$u$ but no nodes of depth $u-e$, but this is harmless.)

Now observe that if $m,n\in\bbz$ with either $m\gs0$ or $n\ls0$, then the reduction of the sequence $[m][n]$ is $[m+n]$. So the reduced $(A,i)$-signature of~$\la$ is the reduction of the sequence
\[
S=\dots\big[\add i{u_3}\la-\rem i{u_3-e}\la\big]\big[\add i{u_2}\la-\rem i{u_2-e}\la\big]\big[\add i{u_1}\la-\rem i{u_1-e}\la\big].
\]

The calculation of the $(B,i)$-signature of~$\mu$ works in exactly the same way, except that we interchange ``above'' and ``below'' when considering nodes of the same depth. So we find that the reduced $(B,i)$-signature of~$\mu$ is the reduction of the sequence
\[
\dots\big[\add i{u_3}\mu-\rem i{u_3-e}\mu\big]\big[\add i{u_2}\mu-\rem i{u_2-e}\mu\big]\big[\add i{u_1}\mu-\rem i{u_1-e}\mu\big].
\]
But \cref{addrem} shows that this sequence coincides with~$S$; moreover, the addable or removable nodes of~$\la$ and~$\mu$ corresponding to each sign in this sequence have the same depth. The result follows.
\end{pf}

Now we are able to prove our main theorem.

\begin{pf}[Proof of \cref{main}]
\cref{bijec} shows that $(ez,yz)$-\regn yields a bijection $\rega\to\regb$, with inverse given by $(ez,yz)$-\deregn. Now take $\la\in\rega$, and let~$\mu$ be the $(ez,yz)$-\regn of~$\la$. Then for $i\in\zez$ \cref{samredsig} shows that $\eea i\la\neq0$ \iff $\eeb i\mu\neq0$. Assuming $\eea i\la\neq0$, \cref{samredsig} also shows that $\eea i\la$ and $\eeb i\mu$ are ladder-equivalent. Since $\eea i\la\in\rega$ and $\eeb i\mu\in\regb$, this means that $\eeb i\mu$ is the $(ez,yz)$-\regn of $\eea i\la$ and $\eea i\la$ is the $(ez,yz)$-\deregn of $\eeb i\mu$. So we have an isomorphism of crystals.
\end{pf}

In view of \cref{main}, we extend our notation: if $y\in[1,e-1]$ is rational, then we define \emph{$(e,y)$-\regn} to mean $(ez,yz)$-\regn, where~$z$ is the denominator of~$y$.

\subsection{Composing crystal isomorphisms}\label{composesec}

In this section we extend \cref{main} to show that for any two arm sequences the crystals~$\rega$ and~$\regb$ are isomorphic. We do this by considering finite subcrystals. Given any $n\in\bbn$, let~$\regan$ denote the induced subgraph of the crystal graph of~$\rega$ comprising only partitions of size at most~$ne$.

\begin{lemma}\label{regan}
The crystal~$\regan$ depends only on $A_1,A_2,\dots,A_n$.
\end{lemma}

\begin{pf}
If~$\la$ is a partition of size at most~$ne$, then the hooks of~$\la$ have length at most~$ne$, so in particular~$\la$ cannot have a hook of length~$et$ for $t>n$. So whether~$\la$ lies in~$\rega$ depends only on $A_1,\dots,A_n$. Furthermore, if $(r,c)$ and $(s,d)$ are addable or removable $i$-nodes of~$\la$, then $|r-s-c+d|\ls ne$, so whether $(r,c)\cls(s,d)$ depends only on $A_1,\dots,A_n$. So if $\la\in\regan$ then the $i$-signature of~$\la$, and hence~$\ee i\la$, depend only on $A_1,\dots,A_n$.
\end{pf}

Given two arm sequences~$A$ and~$B$, our strategy will be to construct an isomorphism $\regan\to\regbn$ for every~$n$. Such isomorphisms will automatically be compatible, thanks to the following result on uniqueness of isomorphisms.

\begin{propn}\label{uniqueiso}
Suppose~$A$ is an arm sequence. Then the only isomorphism from~$\rega$ to~$\rega$ is the identity. For any $n\in\bbn$, the only isomorphism from~$\regan$ to~$\regan$ is the identity.
\end{propn}

\begin{pf}
Suppose~$\alpha$ is an isomorphism from~$\rega$ to~$\rega$ or from~$\regan$ to~$\regan$. We prove that $\alpha(\la)=\la$ for every~$\la$ by induction on~$\card\la$. As mentioned in \cref{cryfamilysec},~$\varnothing$ is the unique source in~$\rega$. Because $\card{\ee i\la}<\card\la$ whenever $\ee i\la\neq0$, this means that~$\varnothing$ is also the unique source in~$\regan$. So~$\varnothing$ must be preserved by~$\alpha$. If $\la\neq\varnothing$, then~$\la$ can be written as $\ff{i_1}\dots\ff{i_r}\varnothing$ for some $i_1,\dots,i_r\in\zez$, which means that
\[
\alpha(\la)=\ff{i_1}\dots\ff{i_r}\alpha(\varnothing)=\ff{i_1}\dots\ff{i_r}\varnothing=\la.\qedhere
\]
\end{pf}

\begin{cory}\label{ison}
Suppose~$A$ and~$B$ are arm sequences, that $m<n$ are natural numbers, and $\alpha_m:\regar m\to\regbr m$ and $\alpha_n:\regan\to\regbn$ are isomorphisms. Then $\alpha_n(\regar m)=\regbr m$, and $\alpha_n|_{\regar m}=\alpha_m$.
\end{cory}

\begin{pf}
Obviously $\alpha_m(\varnothing)=\alpha_n(\varnothing)=\varnothing$, since~$\varnothing$ is the unique source in any of these crystals. This then implies that $\card{\alpha_n(\la)}=\card\la$ for every $\la\in\regan$, since~$\card\la$ is the length of every directed path from~$\varnothing$ to~$\la$ in~$\rega$. Hence $\alpha_n(\regar m)\subseteq\regbr m$; then the same argument with~$\alpha_n^{-1}$ in place of~$\alpha$ shows that $\alpha_n(\regar m)=\regbr m$.

So $\alpha_n|_{\regar m}$ gives an isomorphism from~$\regar m$ to~$\regbr m$. By \cref{uniqueiso} there is only one isomorphism from~$\regar m$ to~$\regar m$, and so there is at most one isomorphism from~$\regar m$ to~$\regbr m$. So $\alpha_n|_{\regar m}=\alpha_m$.
\end{pf}

So if we can construct an isomorphism $\regan\to\regbn$ for every~$n$, then we obtain an isomorphism $\rega\to\regb$ by gluing these isomorphisms together. The next lemma shows how to realise any given~$\regan$ as~$\rgr{A^{y+}}n$ for some $y\in\bbq$.

\begin{lemma}\label{yfromb}
Suppose~$A$ is an arm sequence and $n\in\bbn$, and let
\[
y=\max\lset{\mfrac{A_t}t}{t\in\{1,\dots,n\}}.
\]
Then $A_t=\inp{yt}$ for $t=1,\dots,n$. So $\regan=\rgr{A^{y+}}n$.
\end{lemma}

\begin{pf}
The choice of~$y$ means that $A_t\ls yt$ for each~$t$. So we need to show that $A_t>yt-1$ for each~$t$. If instead $A_t\ls yt-1$ for some~$t$, then for any~$u$ the definition of an arm sequence implies that $A_{tu}\ls A_tu+u-1<ytu$. But if we choose~$u$ so that $A_u=yu$, then the definition of an arm sequence also gives $A_{tu}\gs A_ut=ytu$, a contradiction.

So $(A_1,\dots,A_t)=(A^{y+}_1,\dots,A^{y+}_t)$, and so by \cref{regan} $\regan=\rgr{A^{y+}}n$.
\end{pf}

Now we can prove the main result of this section.

\begin{thm}\label{maincompose}
Suppose~$A$ and~$B$ are arm sequences. Then~$\regan$ and~$\regbn$ are isomorphic for every $n\in\bbn$, and so~$\rega$ is isomorphic to~$\regb$.
\end{thm}

\begin{pf}
Recall the total order on arm sequences defined by $A\gs B$ if $A_t\gs B_t$ for all~$t$. We will assume $A\gs B$, and show that~$\regan$ is isomorphic to~$\regbn$ by induction on the natural number $N=\sum_{t=1}^n(A_t-B_t)$. If $N=0$, then $(A_1,\dots,A_n)=(B_1,\dots,B_n)$, and so the identity map is an isomorphism from~$\regan$ to~$\regbn$. If not, then let $y=\max\lset{\mfrac{A_t}t}{t\in\{1,\dots,n\}}$. Then by \cref{yfromb} $\regan=\rgr{A^{y+}}n$.

Since by assumption $A_t>B_t$ for at least one $t\in\{1,\dots,n\}$, we get $A^{y+}>B$, so that $A^{y-}\gs B$. Furthermore, $A^{y-}_t<A^{y+}_t$ for at least one $t\in\{1,\dots,n\}$ (namely, any value of~$t$ for which $A_t=yt$) so $\sum_{t=1}^n(A^{y-}_t-B_t)<N$, so by induction $\rgr{A^{y-}}n$ is isomorphic to~$\regbn$. By \cref{main}~$\rg{A^{y+}}$ and~$\rg{A^{y-}}$ are isomorphic via a \regn map, and so $\regan=\rgr{A^{y+}}n$ is isomorphic to~$\rgr{A^{y-}}n$, and hence to~$\regbn$.

So~$\regan$ is isomorphic to~$\regbn$ for every~$n$, and so (as explained at the start of this section)~$\rega$ is isomorphic to~$\regb$.
\end{pf}

\begin{eg}
Suppose $e=4$, and~$A$ and~$B$ are arm sequences with $(A_1,A_2,A_3,A_4)=(2,4,6,8)$ and $(B_1,B_2,B_3,B_4)=(1,2,4,5)$. We consider the subcrystals~$\rgr A4$ and~$\rgr B4$ obtained by taking partitions of size at most~$16$. We can show that these are isomorphic via the chain of isomorphisms in \cref{figchain} (where for each equality of crystals we give the common first four entries of the corresponding arm sequences).

Under the isomorphism $\rgr A4\to\rgr B4$, the partition $(4,3^2,2,1^4)$ maps to $(6,4,2,1^4)$, via the sequence also shown in \cref{figchain}.

\begin{figure}[ht]
\[
\begin{tikzpicture}[yscale=-2,xscale=2.6,rotate=90]
\draw(0,5)node(a){$\rgr A4$};
\draw(0,4)node(2p){$\rgr{A^{2+}}4$};
\draw(1,5)node(2m){$\rgr{A^{2-}}4$};
\draw(1,4)node(74p){$\rgr{A^{7/4+}}4$};
\draw(2,5)node(74m){$\rgr{A^{7/4-}}4$};
\draw(2,4)node(53p){$\rgr{A^{5/3+}}4$};
\draw(3,5)node(53m){$\rgr{A^{5/3-}}4$};
\draw(3,4)node(32p){$\rgr{A^{3/2+}}4$};
\draw(4,5)node(32m){$\rgr{A^{3/2-}}4$};
\draw(4,4)node(b){$\rgr B4$};
\draw[double,double distance=2pt](a)--(2p)node[midway,above]{\footnotesize$(2,4,6,8)$};
\draw[double,double distance=2pt](2m)--(74p)node[midway,above]{\footnotesize$(1,3,5,7)$};
\draw[double,double distance=2pt](74m)--(53p)node[midway,above]{\footnotesize$(1,3,5,6)$};
\draw[double,double distance=2pt](53m)--(32p)node[midway,above]{\footnotesize$(1,3,4,6)$};
\draw[double,double distance=2pt](32m)--(b)node[midway,above]{\footnotesize$(1,2,4,5)$};
\draw[->](2p)--(2m)node[midway,right]{\ \footnotesize$(4,2)$-\regn};
\draw[->](74p)--(74m)node[midway,right]{\ \footnotesize$(16,7)$-\regn};
\draw[->](53p)--(53m)node[midway,right]{\ \footnotesize$(12,5)$-\regn};
\draw[->](32p)--(32m)node[midway,right]{\ \footnotesize$(8,3)$-\regn};
\draw(0,2)node(a){$(4,3^2,2,1^4)$};
\draw(1,2)node(b){$(5,4,2,1^5)$};
\draw(2,2)node(c){$(5,4,2,1^5)$};
\draw(3,2)node(d){$(6,4,2,1^4)$};
\draw(4,2)node(e){$(6,4,2,1^4)$};
\draw[->](a)--(b)node[midway,right]{\footnotesize$(4,2)$-\regn};
\draw[->](b)--(c)node[midway,right]{\footnotesize$(16,7)$-\regn};
\draw[->](c)--(d)node[midway,right]{\footnotesize$(12,5)$-\regn};
\draw[->](d)--(e)node[midway,right]{\footnotesize$(8,3)$-\regn};
\end{tikzpicture}
\]
\caption{}\label{figchain}
\end{figure}
\end{eg}

We end this section with some comments on the case $e=2$. In this case there are only two arm sequences
\[
A^{1+}=(1,2,3,\dots),\qquad A^{1-}=(0,1,2,\dots),
\]
so our family of crystals contains nothing beyond the Misra--Miwa model $\rg{A^{1+}}$ and its dual $\rg{A^{1-}}$. (Berg's work remains valid in the case $e=2$, but his model now coincides with the Misra--Miwa crystal.) \cref{main} remains valid when $e=2$: the only possible value of~$y$ now is $y=1$, and the isomorphisms between~$\rg{A^{1+}}$ and~$\rg{A^{1-}}$ are given by $(2,1)$-\regn and -\deregn. In fact if~$\la$ is a $2$-restricted partition then its $(2,1)$-\regn is simply~$\la'$, and similarly for $(2,1)$-\deregn.

So the main results of the present section remain valid for $e=2$, though they tell us almost nothing. But we will need to consider the case $e=2$ later.

\section{The Mullineux map}\label{mullsec}

Now we use the results of the previous section to give a new combinatorial algorithm for the Mullineux map.

\subsection{Background on the Mullineux map}

The Mullineux map first arose in the representation theory of the symmetric group~$\sss n$. When~$e$ is a prime number and~$\bbf$ is a field of characteristic~$e$, the irreducible~$\bbf\sss n$-modules are the \emph{James modules}~$\jms\la$, for the different $e$-regular partitions~$\la$ of~$n$; these were first constructed in \cite{j0}. One of these modules is the one-dimensional \emph{sign module}~$\sgn$, on which a permutation $\pi\in\sss n$ acts as $\sgn(\pi)$. Tensoring an irreducible module with a one-dimensional module yields an irreducible module, which means that there is a function~$\mule$ (the \emph{Mullineux involution}) on the set of $e$-regular partitions of~$n$ such that
\[
\jms\la\otimes\sgn\cong\jms{\mule(\la)}
\]
for each~$\la$. The \emph{Mullineux problem} is to give a combinatorial description of~$\mule$, and there are now several known solutions to this problem. Mullineux \cite{mull} gave a combinatorial map~$\Mule$, and conjectured that $\Mule=\mule$. Meanwhile, Kleshchev \cite{kles} gave a different combinatorial map~$\kule$, and proved (using his modular branching rules for the symmetric groups) that $\kule=\mule$. Shortly afterwards, Ford and Kleshchev \cite{fokl} proved the purely combinatorial result that $\kule=\Mule$, thus proving Mullineux's conjecture. Shorter proofs of this statement was given by Bessenrodt and Olsson \cite{beol} and by Xu \cite{xu}, who gave another combinatorial map~$\xule$ and showed that $\kule=\xule=\Mule$. Finally, Brundan and Kujawa \cite{brku} gave yet another combinatorial map~$\sule$ based on Serganova's work on the general linear supergroup, and showed that $\mule=\sule=\xule$.

None of the combinatorial parts of this story depend on~$e$ being prime, and all the combinatorial maps mentioned in the last paragraph are defined when~$e$ is any integer great then~$1$. In fact there is an algebraic interpretation of this more general situation: the Mullineux involution for arbitrary~$e$ can be defined by generalising the construction above from the group algebra of~$\sss n$ in characteristic~$e$ to an Iwahori--Hecke algebra of~$\sss n$ in quantum characteristic~$e$; Brundan \cite{brun} showed how to generalise Kleshchev's results to this setting, so that this more general Mullineux map is given by all the combinatorial algorithms mentioned in the last paragraph.

\subsection{Crystals and a new algorithm for the Mullineux map}\label{newalgsec}

The Brundan--Kleshchev approach to the Mullineux problem also provides a link with crystals. Define a \emph{signed isomorphism} between crystals~$A$ and~$B$ to be a bijection $\alpha:A\to B$ such that $\ff ib=c$ \iff $\ff{-i}\alpha(b)=\alpha(c)$, for $b,c\in A$ and $i\in\zez$. In other words,~$\alpha$ is an isomorphism of directed graphs under which each arrow $\stackrel i\longrightarrow$ maps to an arrow~$\stackrel{-i}\longrightarrow$.

Recall that~$\blo$ denotes the crystal of the irreducible highest-weight module for $\widehat{\mathfrak{sl}}_e$ with highest weight~$\La_0$. There is a signed isomorphism $\blo\to \blo$; this reflects the diagram automorphism of~$\widehat{\mathfrak{sl}}_e$ preserving the weight~$\La_0$. For each of our combinatorial models~$\rega$ for this crystal, this signed isomorphism is realised as an involutory bijection on the set of $A$-regular partitions. In the particular case of the arm sequence $A=A^{1-}=(0,1,2,\dots)$, the crystal~$\rega$ consists of all $e$-regular partitions, and the Brundan--Kleshchev branching rules show that the signed automorphism in this case (which is unique, in view of \cref{uniqueiso}) is the Mullineux map.

We can use this to realise the Mullineux map in terms of (unsigned) isomorphisms between different crystal models. The following result comes directly from the definitions by conjugating partitions (and using the fact that the residue of the node $(c,r)$ is the negative of the residue of $(r,c)$).

\begin{propn}\label{crystaliso}
Suppose~$A$ is an arm sequence, and define the conjugate arm sequence~$A'$ by
\[
A'_t=et-A_t-1.
\]
If~$\la$ is a partition, then $\la\in\rega$ \iff $\la'\in\rg{A'}$. Furthermore, there is a signed isomorphism $\rega\to\rg{A'}$ given by $\la\mapsto\la'$.
\end{propn}

If we consider in particular the arm sequence $A=A^{1-}$, then $A'=A^{(e-1)+}$, and the $A'$-regular partitions are precisely the $e$-restricted partitions. \cref{crystaliso} yields a signed isomorphism $\rega\to\rg{A'}$. Composing this with the unique isomorphism $\rg{A'}\to\rega$ gives a signed isomorphism from $\rega\to\rega$, which must be the Mullineux map, by uniqueness. So we have proved the following.

\begin{propn}\label{mullconj}
Let~$\alpha$ be the unique isomorphism from~$\rg{A^{(e-1)+}}$ to~$\rg{A^{1-}}$. Then
\[
\mule(\la)=\alpha(\la')
\]
for any \erp~$\la$.
\end{propn}

The results of \cref{composesec} allow us to express the isomorphism~$\alpha$ as a composition of regularisation maps. We use this description to give a new combinatorial algorithm to compute~$\mule(\la)$, for any \erp~$\la$. Recall that if $y\in[1,e-1]\cap\bbq$, then we use the term \earegn to mean $(ez,yz)$-\regn, where~$z$ is the denominator of~$y$. Now we can give our algorithm for computing~$\mule(\la)$.

\begin{enumerate}
\item
Let $\mu=\la'$ and $x=e-1$.
\item
If there exists a rational number $y\in[1,x]$ such that~$\mu$ is not $A^{y-}$-regular, then let~$y$ be the largest such number. Replace~$\mu$ with its \earegn and replace~$x$ with~$y$. 
\item
Repeat Step 2 until~$\mu$ is $A^{y-}$-regular for every $y\in[1,x]\cap\bbq$.
\item
Output~$\mu$.
\end{enumerate}

In other words, we start with the partition $\mu=\la'$, and imagine a variable~$y$ decreasing continuously from $e-1$ to~$1$; for each rational value of~$y$, we replace~$\mu$ with its $(e,y)$-\regn.

In Step 2 of the algorithm, finding the largest $y\in[1,x]$ such that~$\mu$ is not $A^{y-}$-regular is actually straightforward. Suppose~$\mu$ has a hook with length~$re$ and arm length~$a$. The existence of this hook means that~$\mu$ fails to be $A^{y-}$-regular whenever $a/r<y\ls(a+1)/r$. So if we define the \emph{slope} of this hook to be $(a+1)/r$, then the required value of~$y$ is the largest $y\in[1,x]$ which occurs as the slope of an $re$-hook in~$\mu$.

\begin{eg}
Suppose $e=3$ and $\la=(6,2,1)$. We set $\mu=\la'=(3,2,1^4)$ and $x=2$. We draw the Young diagram of~$\mu$; for the hooks of length divisible by~$3$, we fill the corresponding node with the slope of the hook:
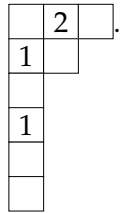
\begin{align*}
&\young(\ 2\ ,1\ ,\ ,1,\ ,\ ).
\\
\intertext{In this case the largest $y\ls2$ for which~$\mu$ is not $A^{y+}$-regular is $y=2$.  So we replace~$\mu$ with its $(3,2)$-\regn $(4,1^5)$, for which we draw a similar diagram:}
&\young({{\tsfrac43}}3\ \ ,\ ,\ ,1,\ ,\ ).
\\
\intertext{The existence of a $9$-hook with arm length~$3$ means that $y=\frac43$. So we replace~$\mu$ with its $(3,\frac43)$-\regn (that is, its $(9,4)$-\regn) $(5,1^4)$, and draw its diagram:}
&\young({{\tsfrac53}}\ 3\ \ ,\ ,1,\ ,\ ).
\\
\intertext{We set $y=1$, and replace~$\mu$ with its $(3,1)$-\regn $(5,2^2)$.}
&\young(\ 23\ \ ,2\ ,\ \ ).
\end{align*}
Since~$y$ is now~$1$, we output $\mul3(6,2,1)=(5,2^2)$.
\end{eg}

\subsection{Regularisation and the Mullineux map}

In this section we use our new algorithm to give a new proof of the main result from the author's paper \cite{mfmull} relating the Mullineux map and $e$-regularisation. Given a partition~$\la$, we write~$\la\reg$ for the $(e,1)$-\regn of~$\la$. Given a hook of~$\la$ with arm length~$a$ and leg length~$l$, we say that the hook is \emph{steep} if $l\gs(e-1)a$, or \emph{shallow} if $a\gs(e-1)l$. Now we have the following result; this generalises a result of Bessenrodt, Olsson and Xu \cite[Theorem 4.8]{box} and  confirms a conjecture of Lyle \cite{lyle}.

\begin{thmc}{mfmull}{Conjecture 1.6}\label{mullregmain}
Suppose~$\la$ is a partition. Then
\[
\mule(\la\reg)\dom(\la')\reg,
\]
with equality \iff every hook of~$\la$ with length divisible by~$e$ is either steep or shallow.
\end{thmc}

\begin{rmks}\indent
\begin{enumerate}
\vspace*{-\topsep}
\item
In fact, the statement given above is stronger than that given in \cite{mfmull}, where dominance is not considered; only the statement about when $\mule(\la\reg)=(\la')\reg$ is proved. However, the method of proof in \cite{mfmull} can easily be adapted to show the dominance part of the result as well.
\item
In \cite{mfmull}, the ``if'' part of \cref{mullregmain} is proved combinatorially, while an algebraic proof (using canonical basis coefficients) is given for the ``only if'' part. The proof we give here is purely combinatorial.
\end{enumerate}
\end{rmks}

First we need a lemma.

\begin{lemma}\label{ladderdom}
Suppose $e-1\gs y>1$. Suppose~$\la$ is an $(e,y)$-singular partition, and let~$\mu$ be its $(e,y)$-\regn. Then $\mu\reg\doms\la\reg$.
\end{lemma}

\begin{pf}
Each $(e,1)$-ladder contains exactly one node in each column. So we can number the $(e,1)$-ladders $\dots,\call_1,\call_2,\dots$ by letting~$\call_r$ be the ladder containing $(r,1)$. We say that ladder~$\call_s$ is \emph{later} than ladder~$\call_r$ if $s>r$. Then each node in ladder~$\call_r$ has a node in~$\call_{r+1}$ immediately below it, and as a consequence the number of nodes in~$\call_{r+1}\setminus\la$ is at least the number of nodes in $\call_r\setminus\la$.

Now take $M>0$ such that $\la,\mu\subseteq\call_1\cup\dots\cup\call_M$, let $N=\card{\call_M\cap\bbn^2}$, and define
\[
\hat\la_r=N-\card{\call_r\setminus\la}
\]
for $r=1,\dots,N$, with $\hat\la_r=0$ for $r>N$. Define~$\hat\mu$ similarly. The previous paragraph shows that~$\hat\la$ and~$\hat\mu$ are partitions, and obviously $\card{\hat\la}=\card{\hat\mu}$.

Now observe that if $(r,c)$ and $(s,d)$ are nodes in the same $(e,y)$-ladder with $s<r$, then $(s,d)$ lies in a later $(e,1)$-ladder than $(r,c)$. Since~$\mu$ is constructed from~$\la$ by moving nodes up their $(e,y)$-ladders, this means that $\hat\la\doms\hat\mu$.

Now consider~$(\la\reg)'$. Since~$\la\reg$ is constructed simply by moving all the nodes of~$\la$ up to the highest positions in their $(e,1)$-ladders,~$(\la\reg)'_d$ equals the number of~$r$ for which $\card{\call_r\setminus\la}<d$. Hence for any $d\gs1$,
\begin{align*}
(\la\reg)'_1+\dots+(\la\reg)'_d&=\sum_{i=1}^d\cardx{\vphantom{a^l}\lset{r\gs1}{\hat\la_r>N-i}}\\
&=\hat\la'_{N-d+1}+\dots+\hat\la'_N\\
&=\card{\hat\la}-(\hat\la'_1+\dots+\hat\la'_{N-d})
\end{align*}
and similarly for~$\mu$. Now the fact that $\hat\la\doms\hat\mu$ gives $(\la\reg)'\doms(\mu\reg)'$, and hence $\la\reg\domsby\mu\reg$.
\end{pf}

\begin{pf}[Proof of \cref{mullregmain}]
Suppose first that every hook of~$\la$ with length divisible by~$e$ is either shallow or steep. This means in particular that~$\la'$ is $(e,e-1)$-regular. Since~$\la$ and~$\la\reg$ lie (by definition) in the same $(e,1)$-ladder class,~$\la'$ and~$(\la\reg)'$ lie in the same $(e,e-1)$-ladder class. So~$\la'$ is the $(e,e-1)$-\regn of~$(\la\reg)'$.

When we apply our algorithm to compute~$\mule(\la\reg)$, we first compute~$(\la\reg)'$ and then its $(e,e-1)$-\regn, which is~$\la'$, and then apply $(e,y)$-\regn for every rational $y<e-1$. But because every hook of~$\la$ with length divisible by~$e$ is either shallow or steep,~$\la'$ is $(e,y)$-regular whenever $e-1>y>1$. So none of the $(e,y)$-\regn maps have any effect until we reach $y=1$, at which point we replace~$\la'$ with its $(e,1)$-\regn~$(\la')\reg$. So $\mule(\la\reg)=(\la')\reg$.

Conversely, suppose~$\la$ has an $er$-hook which is neither steep nor shallow, for some~$r$. Suppose first that~$\la'$ is $(e,e-1)$-regular. Then the calculation of~$\mule(\la\reg)$ begins as in the case above, by computing~$(\la\reg)'$ and then $(e,e-1)$-regularising to get~$\la'$. Then we reach~$\mule(\la\reg)$ by applying $(e,y)$-\regn for some finite list of values $y=y_1,\dots,y_r$, with $e-1>y_1>\dots>y_r>1$, and then finally (if necessary) applying $(e,1)$-regn. So we construct a list of partitions
\[
\la'=\mu(0)\neq\mu(1)\neq\dots\neq\mu(r),
\]
where~$\mu(k)$ is the $(e,y_k)$-\regn of $\mu(k-1)$ for each~$k$, and $\mule(\la\reg)=\mu(r)\reg$.

The fact that~$\la$ (and hence~$\la'$) has an $er$-hook which is neither steep nor shallow means that~$\la'$ is $(e,y)$-singular for some~$y$ with $e-1>y>1$. (Indeed, suppose~$\la'$ has an $et$-hook with arm length~$a$ and leg length~$l$; then~$\la'$ is $(e,y)$-singular for $y=(a-1)/t$, which lies strictly between~$1$ and $e-1$.) Hence $r\gs1$. Now by \cref{ladderdom}
\[
(\la')\reg\domsby\mu(1)\reg\domsby\cdots\domsby\mu(r)\reg=\mule(\la\reg),
\]
as required.

The case where~$\la$ is $(e,e-1)$-singular is similar: in this case we begin by computing~$(\la\reg)'$, and then applying $(e,y)$-\regn for $y=y_1,\dots,y_r$ where now $e-1=y_1>\dots>y_r>1$. Again this gives a sequence of partitions $\mu(1),\dots,\mu(r)$ with
\[
\mu(1)\reg\domsby\cdots\domsby\mu(r)\reg=\mule(\la\reg).
\]
In this case~$\mu(1)$ is the $(e,e-1)$-\regn of~$(\la\reg)'$, which is the same as the $(e,e-1)$-\regn of~$\la'$. So \cref{ladderdom} gives $\mu(1)\reg\doms(\la')\reg$, which is what we need.
\end{pf}

\subsection{The Mullineux map and separated partitions on the abacus}\label{splitsec}

In this section we use our new algorithm for the Mullineux map to prove a reduction theorem which shows that for certain partitions (characterised by a property of their \abds) the effect of the Mullineux map can be calculated using Mullineux maps for smaller values of~$e$. Applying this repeatedly yields a new (and purely combinatorial) proof of Paget's theorem \cite{pag} describing the effect of the Mullineux map on partitions labelling modules in RoCK blocks of symmetric groups and Iwahori--Hecke algebras.

We fix a set $I$ of integers which is a union of congruence classes modulo~$e$, and set $\ol I=\bbz\setminus I$. We define~$c$ to be the number of congruence classes contained in~$I$, and set $\ol c=e-c$. We will assume for now that $2\ls c\ls e-2$, but later we will explain how our results can be extended to the cases where $c=1$ or $e-1$. We fix a large integer~$n$ divisible by~$e$. \Abds with~$e$ runners will always be assumed to have~$n$ beads.

Given a partition~$\la$, we take the $n$-bead \abd with~$e$ runners, and construct a new \abd with~$c$ runners by discarding all the positions not in~$I$; this defines a partition which we write as~$\la_I$. Similarly, discarding the positions in~$I$ gives a~$\ol c$-runner \abd for a partition which we write as~$\la_{\ol I}$.

\begin{eg}
Take $e=5$, $I=(0+5\bbz)\cup(2+5\bbz)$ and $\la=(5,3^2,2,1)$. Then $\la_I=(2)$ and $\la_{\ol I}=(2,1)$. We show the \abds for these partitions; in the examples in this section, when drawing $e$-runner \abds, we will use white beads for the positions in~$I$.
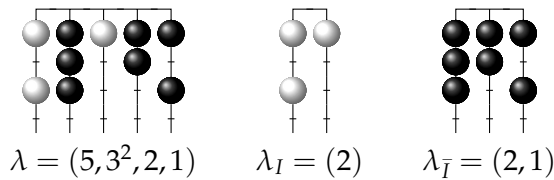
\[
\begin{array}{c@{\qquad}c@{\qquad}c}
\abacus(lmmmr,obobb,nbnbn,obnnb,nnnnn)
&
\abacus(lr,oo,nn,on,nn)
&\abacus(lmr,bbb,bbn,bnb,nnn)
\\
\la=(5,3^2,2,1)&\la_I=(2)&\la_{\ol I}=(2,1)
\end{array}
\]

\end{eg}

Now say that~$\la$ is \emph{\isep} if in the \abd for~$\la$ the first empty position in~$I$ occurs after the last occupied position in~$\ol I$. (This condition is independent of~$n$, given the assumption that~$n$ is divisible by~$e$.) 
Our aim is to describe the effect of the Mullineux map on \isep partitions. In fact we will undertake the equivalent task of computing the composite function $\la\mapsto\mule(\la')$ for $e$-restricted~$\la$; this composite map preserves the $e$-content of a partition, so is better suited to the abacus. So we will take an \isep $e$-restricted partition~$\la$, and describe~$\mule(\la')$ in terms of the partitions~$\la_I$ and~$\la_{\ol I}$ and the maps~$\mul c$ and~$\mul{\ol c}$. Very roughly speaking~$\mule(\la')$ is obtained by replacing~$\la_I$ with~$\mul c(\la_I')$ and~$\la_{\ol I}$ with~$\mul{\ol c}(\la_{\ol I}')$; however, because~$\la_{\ol I}$ need not be~$\ol c$-restricted some additional manipulation is needed.

We will use the \earegn maps discussed earlier, so we start with some simple lemmas looking at the $(e,y)$-regular and $(e,y)$-restricted conditions for \isep partitions.

\begin{lemma}\label{eystuff}
Suppose~$\la$ is an \isep partition.
\begin{enumerate}
\item\label{splitef}
$\la$ is $(e,\ol c)$-regular \iff~$\la_{\ol I}$ is $\ol c$-restricted.
\item\label{splitey}
If~$\la_I$ is $c$-regular, then~$\la$ is $(e,y)$-regular for all $y\in(\ol c,\ol c+1)\cap\bbq$.
\item\label{splitreg}
If $y\in[\ol c+1,e-1]\cap\bbq$, then~$\la$ is $(e,y)$-regular \iff~$\la_I$ is $(c,y-\ol c)$-regular.
\item\label{splitrest}
$\la$ is $e$-restricted \iff~$\la_I$ is $c$-restricted.
\end{enumerate}
\end{lemma}

\begin{pf}
We prove only the ``if'' part of (\ref{splitreg}); the proofs for the other parts are similar.

Suppose~$\la$ is $(e,y)$-singular. Then~$\la$ has a hook with length~$et$ and arm length $yt-1$ for some $t\in\bbn$ for which $yt\in\bbn$. So in the \abd for~$\la$ there is an occupied position $b\gs et$ such that position $b-et$ is empty, and there are exactly~$yt$ empty positions in the range $[b-et,b]$. At most~$\ol ct$ of these positions lie in~$\ol I$, so (because $y>\ol c$) at least one must lie in~$I$. Now the assumption that~$\la$ is \isep means that $b\in I$. Hence the empty position $b-et$ lies in~$I$, so the \isep assumption again means that all the positions in $[b-et,b]\cap\ol I$ are empty. So the $c$-runner \abd for~$\la_I$ has an occupied position~$a$ (corresponding to position~$b$ in the $e$-runner abacus) with position $a-ct$ empty, and exactly $yt-\ol ct$ empty positions in the range $[a-ct,a]$. So~$\nu_I$ is $(c,y-\ol c)$-singular.
\end{pf}

\cref{eystuff}(\ref{splitrest}) allows us to formulate the main theorem of this section, but for this we need some more notation. With our large integer~$n$ fixed as above, we fix an integer $0\ls u\ls n$. We will say that a partition~$\la$ is a \uip if the $n$-bead \abd for~$\la$ has exactly~$u$ beads in positions in~$I$ (and therefore $n-u$ beads in positions in~$\ol I$).  A \uip~$\la$ is determined by the two partitions~$\la_I$ and~$\la_{\ol I}$. Conversely, for any two partitions $\be,\gamma$ there is a unique \uip~$\la$ with $\la_I=\be$ and $\la_{\ol I}=\gamma$.

Given two partitions $\al,\be$, we write $\rgp\al\be$ for the partition
\[
(\ol c\al_1+\be_1,\ol c\al_2+\be_2,\dots)
\]
and $\rsp\al\be$ for the partition obtained by taking the parts of~$\be$ together with~$c$ copies of each of the parts of~$\al$, and arranging these parts in decreasing order.

Now we give the set-up for our main theorem. We keep $I,u$ fixed as above. We take three partitions $\alpha,\be,\gamma$, with~$\be$ being $c$-restricted and~$\gamma$ being $\ol c$-restricted. We define two \uips~$\la$ and~$\mu$ by
\[
\la_I=\be,\qquad\la_{\ol I}=\rgp\alpha\gamma,\qquad\mu_I=\rsp{\alpha'}{\mul c(\be')},\qquad\mu_{\ol I}=\mul{\ol c}(\gamma').
\]
Now we can state the main theorem of this section.

\begin{thm}\label{mainsplit}
Suppose~$\la$ and~$\mu$ defined as above. If~$\la$ and~$\mu$ are both \isep, then $\mule(\la')=\mu$.
\end{thm}

\needspace{5em}
\begin{egs}\indent
\begin{enumerate}
\vspace*{-\topsep}
\item
Take $e=5$ and $I=(1+5\bbz)\cup(4+5\bbz)$, so that $(c,\ol c)=(2,3)$. Taking $n=15$, we choose $u=10$, and
\[
\alpha=(2^2,1),\qquad\be=(2,1^2),\qquad\gamma=(1^3).
\]
We know that $\mul2(3,1)=(3,1)$ and we can calculate $\mul3(3)=(2,1)$, so we get
\[
\la_I=(2,1^2),\qquad\la_{\ol I}=(7^2,4),\qquad\mu_I=(3^3,2^2,1),\qquad\mu_{\ol I}=(2,1),
\]
and hence
\[
\la=(15,11,9,7^3,6,4^3,2,1),\qquad\mu=(17,16,14,10,9,5,2^2,1^2).
\]
We see from the following \abds that~$\la$ and~$\mu$ are both \isep, and we can check that $\mul5(\la')=\mu$.
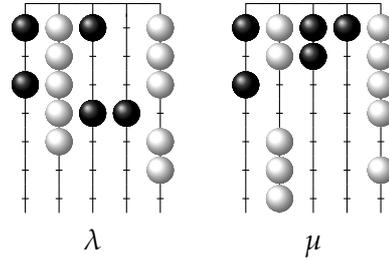
\[
\begin{array}{c@{\qquad}c}
\abacus(lmmmr,bobno,nonno,bonno,nobbn,nonno,nnnno,nnnnn)
&
\abacus(lmmmr,bobbo,nobno,bnnno,nnnno,nonnn,nonno,nonnn)
\\
\la&\mu
\end{array}
\]
\item
We give an example to show that the assumption that~$\la$ and~$\mu$ are both \isep in \cref{mainsplit} is necessary. Take $e=6$, $I=(0+6\bbz)\cup(3+6\bbz)\cup(5+6\bbz)$, $n=12$ and $u=7$. If we take
\[
\alpha=\be=\varnothing,\qquad\gamma=(2^2),
\]
then we get
\[
\la=(2^5),\qquad\mu=(5,2^2,1).
\]
From the \abds we see that~$\la$ is \isep but~$\mu$ is not, and we can check that $\mul6(\la')\neq\mu$.
\[
\begin{array}{c@{\qquad}c}
\abacus(lmmmmr,obbobo,onnobo,obnnnn,nnnnnn)
&
\abacus(lmmmmr,obbobo,obnono,onnnbn,nnnnnn)
\\
\la&\mu
\end{array}
\]
\end{enumerate}
\end{egs}

We will prove \cref{mainsplit} using \cref{mullconj}, which can be rephrased as saying that the map $\la\mapsto\mule(\la')$ is the unique isomorphism between the crystals~$\rg{A^{(e-1)+}}$ and~$\rg{A^{1-}}$. This map is given by applying \earegn for all rational numbers in $[1,e-1]$ in decreasing order. In fact we factorise this isomorphism as a product of three isomorphisms
\[
\rg{A^{(e-1)+}}\longrightarrow\rg{A^{\ol c+}}\longrightarrow\rg{A^{\ol c-}}\longrightarrow\rg{A^{1-}}
\]
and deal with each of these separately. For this we will need to define two intermediate partitions. Keeping the notation from above, we define~$\pi$ and~$\rho$ to be the \uips given by
\[
\pi_I=\mul f(\be'),\qquad\pi_{\ol I}=\rgp\alpha\gamma,\qquad\rho_I=\rsp{\alpha'}{\mul f(\be')},\qquad\rho_{\ol I}=\gamma.
\]
We fix this notation for the rest of this section, and assume from now on that~$\la$ and~$\mu$ are \isep.

The results of \cref{composesec} show that the isomorphism $\rg{A^{(e-1)+}}\longrightarrow\rg{A^{\ol c+}}$ is defined by applying \earegn for all rational $y\in(\ol c,e-1]$ in decreasing order. The next \lcnamecref{regstep1} looks at what happens when we apply one of these regularisation maps.

\begin{propn}\label{regstep1}
Suppose $y\in[\ol c+1,e-1]\cap\bbq$ and~$\nu$ is an \isep partition. Let~$\xi$ be the $(e,y)$-\regn of~$\nu$. Then~$\xi$ is \isep, $\xi_{\ol I}=\nu_{\ol I}$, and~$\xi_I$ is the $(c,y-\ol c)$-\regn of~$\nu_I$.
\end{propn}

\begin{pf}
By \cref{eystuff}(\ref{splitreg})~$\nu$ is $(e,y)$-regular \iff~$\nu_I$ is $(c,y-\ol c)$-regular, and in this case the result is immediate. So assume~$\nu$ is $(e,y)$-singular. We claim that we can find a partition $\kappa\doms\nu$ such that~$\kappa$ and~$\nu$ are $(e,y)$-equivalent,~$\kappa_I$ and~$\nu_I$ are $(c,y-\ol c)$-equivalent and $\kappa_{\ol I}=\nu_{\ol I}$.

The assumption that~$\nu$ is $(e,y)$-singular means that for some $t\in\bbn$ there is an occupied position~$b$ in the \abd for~$\nu$ such that position $b-et$ is empty, and there are exactly~$yt$ empty positions in the range $[b-et,b]$. We take the largest such~$b$, and we will assume that $t=1$. If this is not the case, then we can just replace $e,c,\ol c,y$ with $et,ct,\ol ct,yt$ (keeping the set~$I$ unchanged). Assuming \wolg that~$n$ is divisible by~$et$, all the hypotheses then still apply with these new parameters, and the partition~$\kappa$ obtained will have the desired properties.

So we assume $t=1$, which means in particular that~$y$ is an integer. Now we construct the partition~$\kappa$ as explained in \cref{regnabsec}. Recall that for this we define~$E$ to be the union of the congruence classes containing the empty positions in the range $[b-e,b]$, and let $b_1<\dots<b_m$ be the occupied positions in $\lset{c\in E}{c\gs b}$. We also let $t_1<t_2<\dots$ be the empty positions in $\lset{c\in\bbn_0\setminus E}{c\gs b}$; then we move the beads at positions~$b_i$ and $t_i-e$ to positions $b_i-e$ and~$t_i$ for $i=1,\dots,d$, where~$d$ is minimal such that $t_d<b_{d+1}$. This then gives a partition $\kappa\doms\nu$ which is $(e,y)$-equivalent to~$\nu$.

The fact that~$\nu$ is \isep means (as in the proof of \cref{eystuff}) that $b\in I$, so that all positions in $[b-e,b]\cap\ol I$ are empty, and hence $\ol I\subseteq E$.  Moreover, $b_1,\dots,b_m\in I\cap E$. Hence when we construct~$\kappa$ from~$\nu$ we make no changes to any positions in~$\ol I$, so $\kappa_{\ol I}=\nu_{\ol I}$, as required. Furthermore, the changes we make in the positions in~$I$ correspond exactly to applying the same algorithm (with $e,y$ replaced by $c,y-\ol c$) to~$\nu_I$. Hence~$\kappa_I$ is $(c,y-\ol c)$-equivalent to~$\nu_I$.

So our claim is proved. The fact that $\kappa_I\doms\nu_I$ while $\kappa_{\ol I}=\nu_{\ol I}$ means that~$\kappa$ is also \isep: the last occupied position in the abacus in~$\ol I$ is the same for~$\kappa$ as for~$\nu$, while the fact that $\kappa\doms\nu$ implies that $\kappa'_1\ls\nu'_1$, so the first empty position in~$I$ in the \abd for~$\kappa$ is the same as, or later than, the first empty position in~$I$ in the \abd for~$\nu$.

So the hypotheses of the \lcnamecref{regstep1} apply with~$\nu$ replaced by~$\kappa$. By induction on the dominance order the \lcnamecref{regstep1} holds for~$\kappa$, and so it holds for~$\nu$.
\end{pf}

Applying \cref{regstep1} repeatedly allows us to complete the first step in the proof of \cref{mainsplit}. We recall the partitions $\la,\pi$ from above.

\begin{cory}\label{firstiso}
Let~$\phi$ be the unique isomorphism from~$\rg{A^{(e-1)+}}$ to~$\rg{A^{(\ol c+1)-}}$, and~$\psi$ the unique isomorphism from~$\rg{A^{(e-1)+}}$ to~$\rg{A^{\ol c+}}$. Then $\phi(\la)=\psi(\la)=\pi$.
\end{cory}

\begin{pf}
Using the results in \cref{composesec}, the isomorphism~$\phi$ is given by applying \earegn for all numbers $y\in[\ol c+1,e-1]\cap\bbq$ in decreasing order. Applying \cref{regstep1} at every step, we find that $\phi(\la)_{\ol I}=\la_{\ol I}$, while~$\phi(\la)_I$ is obtained from~$\la_I$ by applying $(c,y)$-\regn for all $y\in[1,c-1]\cap\bbq$ in decreasing order. So (from \cref{newalgsec}) $\phi(\la)_I=\mul c(\la_I')$. Finally,~$\la$ and~$\phi(\la)$ have the same $e$-content, so by \cref{contab}~$\phi(\la)$ is a \uip. So $\phi(\la)=\pi$.

The isomorphism $\psi\circ\phi^{-1}$ from~$\rg{A^{(\ol c+1)-}}$ to~$\rg{A^{\ol c+}}$ is given by applying \earegn for all numbers $y\in(\ol c,\ol c+1)\cap\bbq$ in decreasing order. But~$\pi_I$ is $c$-regular, so~$\pi$ is $(e,y)$-regular for all $y\in(\ol c,\ol c+1)\cap\bbq$, by \cref{eystuff}(\ref{splitey}). So the isomorphism $\psi\circ\phi^{-1}$ fixes~$\pi$, and therefore $\psi(\la)=\pi$.
\end{pf}

Reversing the roles of~$\la$ and~$\mu'$ and applying conjugation throughout, we can complete the third step in the factorisation as well.

\begin{cory}\label{thirdstep}
Let~$\chi$ be the unique isomorphism from~$\rg{A^{\ol c-}}$ to~$\rg{A^{1-}}$. Then $\chi(\rho)=\mu$.
\end{cory}

It remains to show that the isomorphism $\rg{A^{\ol c+}}\to\rg{A^{\ol c-}}$ sends~$\pi$ to~$\rho$. This isomorphism is given by $(e,\ol c)$-\regn. For this we use the following \lcnamecref{efstep}.

\begin{propn}\label{efstep}
Suppose~$\nu$ is an $I$-separated \uip. Suppose~$\nu_{\ol I}$ is not $\ol c$-restricted, and let $s\gs1$ be minimal such that $(\nu_{\ol I})_s-(\nu_{\ol I})_{s+1}\gs\ol c$. Write $\nu_{\ol I}=\rgp{(1^s)}\tau$, and define a \uip~$\xi$ by
\[
\xi_I=\rsp{(s)}{\nu_I},\qquad\xi_{\ol I}=\tau.
\]
Then~$\nu$ and~$\xi$ are $(e,\ol c)$-equivalent.
\end{propn}

\begin{pf}
By \cref{eystuff}(\ref{splitef}), the fact that~$\nu_{\ol I}$ is not $e$-restricted means that~$\nu$ is not $(e,\ol c)$-regular. We apply a single step of the algorithm in \cref{regnabsec} to construct the $(e,\ol c)$-\regn of~$\nu$ on the abacus; we will show that this step produces the partition~$\xi$.

We can find a position~$b$ in the \abd for~$\nu$ and $t\in\bbn$ such that position $b-te$ is empty and there are~$\ol ct$ empty positions in the range $[b-et,b]$. We take the maximal such~$b$. If $b\in I$, then the \isep condition means that all the positions in $[b-et,b]\cap \ol I$ are empty, so that there are at least $\ol ct+1$ empty positions in $[b-et,b]$, a contradiction. So $b\in\ol I$. This means that all the positions in $[b-et,b]\cap I$ are occupied, and therefore all the positions in $[b-et,b-1]\cap\ol I$ are empty. In particular, we can assume $t=1$, and the set~$E$ defined in \cref{regnabsec} now coincides with~$\ol I$.

Consider the $\ol c$-runner \abd~$D$ obtained by discarding all the positions in~$I$, and let~$a$ be the position in~$D$ corresponding to position~$b$ in the $e$-runner \abd for~$\nu$. From the way \abds are constructed, $a=(\nu_{\ol I})_r+n-u-r$ for some~$r$. All the positions in $[a-\ol c,a-1]$ are empty in~$D$, which means that $(\nu_{\ol I})_r-(\nu_{\ol I})_{r+1}\gs\ol c$. The choice of~$b$ means that~$r$ is minimal with this property, so $r=s$. So there are~$s$ occupied positions at or after position~$a$ in~$D$, corresponding to the first~$s$ parts of~$\nu_{\ol I}$. If we let $b_1<\dots<b_s$ be the corresponding positions in the \abd for~$\nu$, then $b_1,\dots,b_s$ are precisely the occupied positions in~$E$ at or after position~$b$. Now let $t_1<t_2<\dots$ be the empty positions in~$I$ after position~$b$. The \isep condition means that $b_s<t_1$, so (according to the recipe in \cref{regnabsec}) we define a new partition~$\kappa$ by moving a bead from position~$b_i$ to position $b_i-e$ and moving a bead from position $t_i-e$ to position~$t_i$, for $i=1,\dots,s$ in turn. The effect of moving the bead from~$b_i$ to $b_i-e$ for each~$i$ is to reduce each of the first~$s$ parts of~$\nu_I$ by~$\ol c$; so $\kappa_{\ol I}=\tau$. By \cref{ablem} the effect of moving a bead from $t_i-e$ to~$t_i$ for each~$i$ is to increase the first~$s$ columns of~$\nu_I$ by~$c$. So $\kappa_I=\rsp{(s)}{\nu_I}=\xi_I$, and therefore $\kappa=\xi$. From \cref{regnabsec}~$\kappa$ is $(e,\ol c)$-equivalent to~$\nu$.
\end{pf}

\begin{cory}\label{secondstep}
$\pi$ and~$\rho$ are $(e,\ol c)$-equivalent, and hence~$\rho$ is the $(e,\ol c)$-\regn of~$\pi$.
\end{cory}

\begin{pf}
Let $a=\alpha_1$, and for each $r=0,\dots,a$ define~$\hat\alpha(r)$ to be the partition obtained by removing the last~$r$ non-zero columns from~$\alpha$, and define~$\check\alpha(r)$ to be the partition obtained by removing the first $a-r$ columns from~$\alpha$. Now define \uips $\nu(0),\dots,\nu(a)$ by
\[
\nu(r)_I=\rsp{\check\alpha(r)'}{\pi_I},\qquad\nu(r)_{\ol I}=\rgp{\hat\alpha(r)}\gamma.
\]
Then $\nu(0)=\pi$ and $\nu(a)=\rho$, and by \cref{efstep} $\nu(r-1)$ is $(e,\ol c)$-equivalent to~$\nu(r)$ for $r=1,\dots,a$. So~$\pi$ and~$\rho$ are $(e,\ol c)$-equivalent. By \cref{eystuff}(\ref{splitef})~$\rho$ is $(e,\ol c)$-regular, and so~$\rho$ is the $(e,\ol c)$-\regn of~$\pi$.
\end{pf}

Now we have proved our main theorem: by \cref{firstiso} the isomorphism from~$\rg{A^{(e-1)+}}$ to~$\rg{A^{\ol c+}}$ sends~$\la$ to~$\pi$; by \cref{secondstep} the isomorphism from~$\rg{A^{\ol c+}}$ to~$\rg{A^{\ol c-}}$ sends~$\pi$ to~$\rho$, and by \cref{thirdstep} the isomorphism from~$\rg{A^{\ol c-}}$ to~$\rg{A^{1-}}$ sends~$\rho$ to~$\mu$. Hence the isomorphism from~$\rg{A^{(e-1)+}}$ to~$\rg{A^{1-}}$ sends~$\la$ to~$\mu$, which is the same as saying that $\mule(\la')=\mu$.

\medskip
Now we consider the case where~$c$ or~$\ol c$ equals~$1$. In this case the construction of the partitions in \cref{mainsplit} still makes sense, but we need to clarify our definitions. The only $1$-regular partition is the empty partition~$\varnothing$ (and this is also the only $1$-restricted partition), and we can define the Mullineux map~$\mul1$ by $\mul1(\varnothing)=\varnothing$. With this convention, the construction of the partitions~$\la$ and~$\mu$ makes sense, and in fact \cref{mainsplit} still holds. The proof just needs to be simplified slightly: if $c=1$ then $\la=\pi$ and the second equality in \cref{firstiso} holds trivially; similarly if $\ol c=1$ then $\rho=\mu$ and \cref{thirdstep} is trivial. \cref{secondstep} needs no modification at all, and the steps still combine to give the main result.

\medskip
We end this section by explaining briefly how \cref{mainsplit} can be use to give a new proof of Paget's theorem \cite[Theorem 2.1]{pag} describing the effect of the Mullineux map for a certain class of partitions; these partitions label modules in \emph{RoCK blocks} of symmetric groups and Iwahori--Hecke algebras. (RoCK blocks -- also called \emph{Rouquier blocks} -- are blocks with especially nice properties, which have played an important role in resolving several conjectures for symmetric groups. These blocks are discussed in detail in \cite{jlm}.) In fact the version we give here of Paget's theorem will be slightly more general, but it is not hard to show that it is implied by Paget's version.

We take a partition~$\la$ and a large integer~$n$ divisible by~$e$, and construct the $n$-bead \abd for~$\la$ with~$e$ runners. For each $i\in\zez$ we regard runner~$i$ of the abacus as a $1$-runner \abd, and let~$\la^{(i)}$ be the corresponding partition. The $e$-tuple $\rtup{\la^{(i)}}{i\in\zez}$ is called the \emph{$e$-quotient} of~$\la$. We also let~$u_i$ be the number of beads on runner~$i$ of the abacus. Now we define a total order~$\ls$ on~$\zez$ by setting $i<j$ if either $u_i<u_j$ or $u_i=u_j$ and runner~$i$ lies to the left of runner~$j$. Define a bijection $\sigma:\{1,\dots,e\}\to\zez$ by $\si(1)<\dots<\si(e)$. Now we say that~$\la$ is \emph{$e$-quotient separated} if for every $1\ls k<l\ls e$, the last occupied position on runner~$\si(k)$ in the \abd for~$\la$ is earlier than the first empty position on runner~$\si(l)$. (The definition of $e$-quotient separated partitions was introduced in the case $e=2$ by James and Mathas \cite{jmq-1}, and the general version has been considered several times; for example in \cite[Section 5.4]{mfaltred}.)

If~$\la$ is $e$-quotient separated, then it is easy to check (and in fact is a special case of \cref{eystuff}(\ref{splitrest}) above) that~$\la$ is $e$-restricted \iff $\la^{(\si(e))}=\varnothing$. Now we can state a version of Paget's theorem as follows.

\begin{thm}\label{pagetthm}
Suppose~$\la$ is an $e$-restricted $e$-quotient separated partition. Let~$\mu$ be the partition with the same $e$-core as~$\la$ and with $e$-quotient $\rtup{\mu^{(i)}}{i\in\zez}$ defined by
\[
\mu^{(\si(1))}=\varnothing,\qquad\mu^{(\si(k))}=(\la^{(\si(k-1))})'\text{ for }k=2,\dots,e.
\]
If~$\mu$ is also $e$-quotient separated, then $\mule(\la')=\mu$.
\end{thm}

\begin{eg}
Take $e=4$, and $\la=(11,10,9,8,7,5^2,4,3,2,1^5)$.
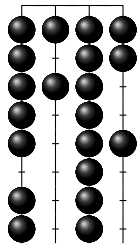
\[
\abacus(lmmr,bbbb,bnbb,bbbn,bnbn,bnbb,nnbn,bnbn,bnbn)
\]
From the \abd we see that~$\la$ is $4$-quotient separated, with
\[
(\si(1),\si(2),\si(3),\si(4))=(1+4\bbz,3+4\bbz,0+4\bbz,2+4\bbz),
\]
\[
(\la^{(\si(1))},\la^{(\si(2))},\la^{(\si(3))},\la^{(\si(4))})=((1),(2),(1^2),\varnothing).
\]
So we define~$\mu$ by
\[
(\mu^{(\si(1))},\mu^{(\si(2))},\mu^{(\si(3))},\mu^{(\si(4))})=(\varnothing,(1),(1^2),(2)),
\]
giving $\mu=(19,10,9,8,7,4,3^3,2,1)$.
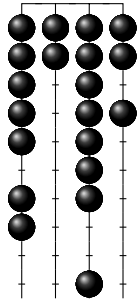
\[
\abacus(lmmr,bbbb,bbbb,bnbn,bnbb,bnbn,nnbn,bnbn,bnnn,nnnn,nnbn)
\]
The \abd shows that~$\mu$ is also $4$-quotient separated, and we can check that $\mul4(\la')=\mu$.
\end{eg}

We can deduce \cref{pagetthm} from \cref{mainsplit} by induction on~$e$. The case $e=2$ is straightforward because~$\mul2$ is the identity map. Assuming $e\gs3$, we apply \cref{mainsplit} with $I=\si(e)$ (so that $c=1$). We take $\beta=\varnothing$ and $\alpha=\la^{(\si(e-1))}$, and we let~$\gamma$ be the partition obtained from~$\la$ by moving all the beads on runner $\si(e-1)$ up to the highest positions on the runner and then deleting runner~$\si(e)$. Then (with $u=u_{\si(e)}$)~$\la$ is the partition defined from~$\alpha$,~$\beta$ and~$\gamma$ in \cref{mainsplit}. Furthermore,~$\gamma$ is $(e-1)$-quotient separated, so we can compute~$\mul{e-1}(\gamma')$ by induction: we find that $\mul{e-1}(\gamma')=\mu_{\ol I}$. Moreover, $\mu_I=\mu^{(e)}=\alpha'$. Hence~$\mu$ is the partition defined from~$\alpha$,~$\beta$ and~$\gamma$ in \cref{mainsplit}, and so the result follows.

\end{document}